\documentclass[reqno,11pt]{amsart}
\newtheorem{theorem}{Theorem}
\newtheorem{lemma}[theorem]{Lemma} 

\newtheorem*{cor}{Corollary}
\DeclareMathOperator{\h}{H}
\begin{document}
\author{Theresia Eisenk\"olbl}

\address{Institut f\"ur Mathematik der Universit\"at Wien\\
Strudlhofgasse 4, A-1090 Wien, Austria.}

\email{teisenko@radon.mat.univie.ac.at}

\title{Rhombus Tilings of a Hexagon with Two Triangles Missing on the
Symmetry Axis}
\begin{abstract}
{We compute the number of rhombus tilings of a 
hexagon with sides $n$, $n$, $N$, $n$, $n$, $N$, where two
triangles on the 
symmetry axis touching in one
vertex are removed. The case of the common vertex being the center of
the hexagon solves a problem posed by Propp.
}
\end{abstract}

\maketitle


\newfont{\fourteenpoint}{cmr10 scaled\magstep3}
\newfont{\fourteenit}{cmti10 scaled\magstep2}
\newfont{\fourteensl}{cmsl10 scaled\magstep2}
\newfont{\fourteensmc}{cmcsc10 scaled\magstep2}
\newfont{\fourteentt}{cmtt10 scaled\magstep2}
\newfont{\fourteenbf}{cmbx10 scaled\magstep2}
\newfont{\fourteeni}{cmmi10 scaled\magstep2}
\newfont{\fourteensy}{cmsy10 scaled\magstep2}
\newfont{\fourteenex}{cmex10 scaled\magstep2}
\newfont{\fourteenmsa}{msam10 scaled\magstep2}
\newfont{\fourteeneufm}{eufm10 scaled\magstep2}
\newfont{\fourteenmsb}{msbm10 scaled\magstep2}

\catcode`\@=11
\font\tenln    = line10
\font\tenlnw   = linew10

\thinlines
\newskip\Einheit \Einheit=0.6cm
\newcount\xcoord \newcount\ycoord
\newdimen\xdim \newdimen\ydim \newdimen\PfadD@cke \newdimen\Pfadd@cke
\PfadD@cke1pt \Pfadd@cke0.5pt
\def\PfadDicke#1{\PfadD@cke#1 \divide\PfadD@cke by2 \Pfadd@cke\PfadD@cke \multiply\PfadD@cke by2}
\long\def\LOOP#1\REPEAT{\def\BODY{#1}\ITERATE}
\def\ITERATE{\BODY \let\next\ITERATE \else\let\next\relax\fi \next}
\let\REPEAT=\fi
\def\Punkt{\hbox{\raise-2pt\hbox to0pt{\hss\scriptsize$\bullet$\hss}}}
\def\DuennPunkt(#1,#2){\unskip
  \raise#2 \Einheit\hbox to0pt{\hskip#1 \Einheit
          \raise-2.5pt\hbox to0pt{\hss\normalsize$\bullet$\hss}\hss}}
\def\NormalPunkt(#1,#2){\unskip
  \raise#2 \Einheit\hbox to0pt{\hskip#1 \Einheit
          \raise-3pt\hbox to0pt{\hss\large$\bullet$\hss}\hss}}
\def\DickPunkt(#1,#2){\unskip
  \raise#2 \Einheit\hbox to0pt{\hskip#1 \Einheit
          \raise-4pt\hbox to0pt{\hss\Large$\bullet$\hss}\hss}}
\def\Kreis(#1,#2){\unskip
  \raise#2 \Einheit\hbox to0pt{\hskip#1 \Einheit
          \raise-4pt\hbox to0pt{\hss\Large$\circ$\hss}\hss}}
\def\Diagonale(#1,#2)#3{\unskip\leavevmode
  \xcoord#1\relax \ycoord#2\relax
      \raise\ycoord \Einheit\hbox to0pt{\hskip\xcoord \Einheit
         \unitlength\Einheit
         \line(1,1){#3}\hss}}
\def\AntiDiagonale(#1,#2)#3{\unskip\leavevmode
  \xcoord#1\relax \ycoord#2\relax \advance\xcoord by -0.05\relax
      \raise\ycoord \Einheit\hbox to0pt{\hskip\xcoord \Einheit
         \unitlength\Einheit
         \line(1,-1){#3}\hss}}
\def\Pfad(#1,#2),#3\endPfad{\unskip\leavevmode
  \xcoord#1 \ycoord#2 \thicklines\ZeichnePfad#3\endPfad\thinlines}
\def\ZeichnePfad#1{\ifx#1\endPfad\let\next\relax
  \else\let\next\ZeichnePfad
    \ifnum#1=1
      \raise\ycoord \Einheit\hbox to0pt{\hskip\xcoord \Einheit
         \vrule height\Pfadd@cke width1 \Einheit depth\Pfadd@cke\hss}%
      \advance\xcoord by 1
    \else\ifnum#1=2
      \raise\ycoord \Einheit\hbox to0pt{\hskip\xcoord \Einheit
        \hbox{\hskip-1pt\vrule height1 \Einheit width\PfadD@cke depth0pt}\hss}%
      \advance\ycoord by 1
    \else\ifnum#1=3
      \raise\ycoord \Einheit\hbox to0pt{\hskip\xcoord \Einheit
         \unitlength\Einheit
         \line(1,1){1}\hss}
      \advance\xcoord by 1
      \advance\ycoord by 1
    \else\ifnum#1=4
      \raise\ycoord \Einheit\hbox to0pt{\hskip\xcoord \Einheit
         \unitlength\Einheit
         \line(1,-1){1}\hss}
      \advance\xcoord by 1
      \advance\ycoord by -1
    \fi\fi\fi\fi
  \fi\next}
\def\hSSchritt{\leavevmode\raise-.4pt\hbox to0pt{\hss.\hss}\hskip.2\Einheit
  \raise-.4pt\hbox to0pt{\hss.\hss}\hskip.2\Einheit
  \raise-.4pt\hbox to0pt{\hss.\hss}\hskip.2\Einheit
  \raise-.4pt\hbox to0pt{\hss.\hss}\hskip.2\Einheit
  \raise-.4pt\hbox to0pt{\hss.\hss}\hskip.2\Einheit}
\def\vSSchritt{\vbox{\baselineskip.2\Einheit\lineskiplimit0pt
\hbox{.}\hbox{.}\hbox{.}\hbox{.}\hbox{.}}}
\def\DSSchritt{\leavevmode\raise-.4pt\hbox to0pt{%
  \hbox to0pt{\hss.\hss}\hskip.2\Einheit
  \raise.2\Einheit\hbox to0pt{\hss.\hss}\hskip.2\Einheit
  \raise.4\Einheit\hbox to0pt{\hss.\hss}\hskip.2\Einheit
  \raise.6\Einheit\hbox to0pt{\hss.\hss}\hskip.2\Einheit
  \raise.8\Einheit\hbox to0pt{\hss.\hss}\hss}}
\def\dSSchritt{\leavevmode\raise-.4pt\hbox to0pt{%
  \hbox to0pt{\hss.\hss}\hskip.2\Einheit
  \raise-.2\Einheit\hbox to0pt{\hss.\hss}\hskip.2\Einheit
  \raise-.4\Einheit\hbox to0pt{\hss.\hss}\hskip.2\Einheit
  \raise-.6\Einheit\hbox to0pt{\hss.\hss}\hskip.2\Einheit
  \raise-.8\Einheit\hbox to0pt{\hss.\hss}\hss}}
\def\SPfad(#1,#2),#3\endSPfad{\unskip\leavevmode
  \xcoord#1 \ycoord#2 \ZeichneSPfad#3\endSPfad}
\def\ZeichneSPfad#1{\ifx#1\endSPfad\let\next\relax
  \else\let\next\ZeichneSPfad
    \ifnum#1=1
      \raise\ycoord \Einheit\hbox to0pt{\hskip\xcoord \Einheit
         \hSSchritt\hss}%
      \advance\xcoord by 1
    \else\ifnum#1=2
      \raise\ycoord \Einheit\hbox to0pt{\hskip\xcoord \Einheit
        \hbox{\hskip-2pt \vSSchritt}\hss}%
      \advance\ycoord by 1
    \else\ifnum#1=3
      \raise\ycoord \Einheit\hbox to0pt{\hskip\xcoord \Einheit
         \DSSchritt\hss}
      \advance\xcoord by 1
      \advance\ycoord by 1
    \else\ifnum#1=4
      \raise\ycoord \Einheit\hbox to0pt{\hskip\xcoord \Einheit
         \dSSchritt\hss}
      \advance\xcoord by 1
      \advance\ycoord by -1
    \fi\fi\fi\fi
  \fi\next}
\def\Koordinatenachsen(#1,#2){\unskip
 \hbox to0pt{\hskip-.5pt\vrule height#2 \Einheit width.5pt depth1 \Einheit}%
 \hbox to0pt{\hskip-1 \Einheit \xcoord#1 \advance\xcoord by1
    \vrule height0.25pt width\xcoord \Einheit depth0.25pt\hss}}
\def\Koordinatenachsen(#1,#2)(#3,#4){\unskip
 \hbox to0pt{\hskip-.5pt \ycoord-#4 \advance\ycoord by1
    \vrule height#2 \Einheit width.5pt depth\ycoord \Einheit}%
 \hbox to0pt{\hskip-1 \Einheit \hskip#3\Einheit 
    \xcoord#1 \advance\xcoord by1 \advance\xcoord by-#3 
    \vrule height0.25pt width\xcoord \Einheit depth0.25pt\hss}}
\def\Gitter(#1,#2){\unskip \xcoord0 \ycoord0 \leavevmode
  \LOOP\ifnum\ycoord<#2
    \loop\ifnum\xcoord<#1
      \raise\ycoord \Einheit\hbox to0pt{\hskip\xcoord \Einheit\Punkt\hss}%
      \advance\xcoord by1
    \repeat
    \xcoord0
    \advance\ycoord by1
  \REPEAT}
\def\Gitter(#1,#2)(#3,#4){\unskip \xcoord#3 \ycoord#4 \leavevmode
  \LOOP\ifnum\ycoord<#2
    \loop\ifnum\xcoord<#1
      \raise\ycoord \Einheit\hbox to0pt{\hskip\xcoord \Einheit\Punkt\hss}%
      \advance\xcoord by1
    \repeat
    \xcoord#3
    \advance\ycoord by1
  \REPEAT}
\def\Label#1#2(#3,#4){\unskip \xdim#3 \Einheit \ydim#4 \Einheit
  \def\lo{\advance\xdim by-.5 \Einheit \advance\ydim by.5 \Einheit}%
  \def\llo{\advance\xdim by-.25cm \advance\ydim by.5 \Einheit}%
  \def\loo{\advance\xdim by-.5 \Einheit \advance\ydim by.25cm}%
  \def\o{\advance\ydim by.25cm}%
  \def\ro{\advance\xdim by.5 \Einheit \advance\ydim by.5 \Einheit}%
  \def\rro{\advance\xdim by.25cm \advance\ydim by.5 \Einheit}%
  \def\roo{\advance\xdim by.5 \Einheit \advance\ydim by.25cm}%
  \def\l{\advance\xdim by-.30cm}%
  \def\r{\advance\xdim by.30cm}%
  \def\lu{\advance\xdim by-.5 \Einheit \advance\ydim by-.6 \Einheit}%
  \def\llu{\advance\xdim by-.25cm \advance\ydim by-.6 \Einheit}%
  \def\luu{\advance\xdim by-.5 \Einheit \advance\ydim by-.30cm}%
  \def\u{\advance\ydim by-.30cm}%
  \def\ru{\advance\xdim by.5 \Einheit \advance\ydim by-.6 \Einheit}%
  \def\rru{\advance\xdim by.25cm \advance\ydim by-.6 \Einheit}%
  \def\ruu{\advance\xdim by.5 \Einheit \advance\ydim by-.30cm}%
  #1\raise\ydim\hbox to0pt{\hskip\xdim
     \vbox to0pt{\vss\hbox to0pt{\hss$#2$\hss}\vss}\hss}%
}
\catcode`\@=12






\def\setRevDate $#1 #2 #3${#2}
\def\TeXdrawId{\setRevDate $Date: 1995/12/19 16:40:42 $ TeXdraw V2R0}
\chardef\catamp=\the\catcode`\@
\catcode`\@=11
\long
\def\centertexdraw #1{\hbox to \hsize{\hss
\btexdraw #1\etexdraw
\hss}}
\def\btexdraw {\x@pix=0             \y@pix=0
\x@segoffpix=\x@pix  \y@segoffpix=\y@pix
\t@exdrawdef
\setbox\t@xdbox=\vbox\bgroup\offinterlineskip
\global\d@bs=0
\global\t@extonlytrue
\global\p@osinitfalse
\s@avemove \x@pix \y@pix
\m@pendingfalse
\global\p@osinitfalse
\p@athfalse
\the\everytexdraw}
\def\etexdraw {\ift@extonly \else
\t@drclose
\fi
\egroup
\ifdim \wd\t@xdbox>0pt
\t@xderror {TeXdraw box non-zero size,
possible extraneous text}%
\fi
\vbox {\offinterlineskip
\pixtobp \xminpix \l@lxbp  \pixtobp \yminpix \l@lybp
\pixtobp \xmaxpix \u@rxbp  \pixtobp \ymaxpix \u@rybp
\hbox{\t@xdinclude 
[{\l@lxbp},{\l@lybp}][{\u@rxbp},{\u@rybp}]{\p@sfile}}%
\pixtodim \xminpix \t@xpos  \pixtodim \yminpix \t@ypos
\kern \t@ypos
\hbox {\kern -\t@xpos
\box\t@xdbox
\kern \t@xpos}%
\kern -\t@ypos\relax}}
\def\drawdim #1 {\def\d@dim{#1\relax}}
\def\setunitscale #1 {\edef\u@nitsc{#1}%
\realmult \u@nitsc \s@egsc \d@sc}
\def\relunitscale #1 {\realmult {#1}\u@nitsc \u@nitsc
\realmult \u@nitsc \s@egsc \d@sc}
\def\setsegscale #1 {\edef\s@egsc {#1}%
\realmult \u@nitsc \s@egsc \d@sc}
\def\relsegscale #1 {\realmult {#1}\s@egsc \s@egsc
\realmult \u@nitsc \s@egsc \d@sc}
\def\bsegment {\ifp@ath
\f@lushbs
\f@lushmove
\fi
\begingroup
\x@segoffpix=\x@pix
\y@segoffpix=\y@pix
\setsegscale 1
\global\advance \d@bs by 1\relax}
\def\esegment {\endgroup
\ifnum \d@bs=0
\writetx {es}%
\else
\global\advance \d@bs by -1
\fi}
\def\savecurrpos (#1 #2){\getsympos (#1 #2)\a@rgx\a@rgy
\s@etcsn \a@rgx {\the\x@pix}%
\s@etcsn \a@rgy {\the\y@pix}}
\def\savepos (#1 #2)(#3 #4){\getpos (#1 #2)\a@rgx\a@rgy
\coordtopix \a@rgx \t@pixa
\advance \t@pixa by \x@segoffpix
\coordtopix \a@rgy \t@pixb
\advance \t@pixb by \y@segoffpix
\getsympos (#3 #4)\a@rgx\a@rgy
\s@etcsn \a@rgx {\the\t@pixa}%
\s@etcsn \a@rgy {\the\t@pixb}}
\def\linewd #1 {\coordtopix {#1}\t@pixa
\f@lushbs
\writetx {\the\t@pixa\space sl}}
\def\setgray #1 {\f@lushbs
\writetx {#1 sg}}
\def\lpatt (#1){\listtopix (#1)\p@ixlist
\f@lushbs
\writetx {[\p@ixlist] sd}}
\def\lvec (#1 #2){\getpos (#1 #2)\a@rgx\a@rgy
\s@etpospix \a@rgx \a@rgy
\writeps {\the\x@pix\space \the\y@pix\space lv}}
\def\rlvec (#1 #2){\getpos (#1 #2)\a@rgx\a@rgy
\r@elpospix \a@rgx \a@rgy
\writeps {\the\x@pix\space \the\y@pix\space lv}}
\def\move (#1 #2){\getpos (#1 #2)\a@rgx\a@rgy
\s@etpospix \a@rgx \a@rgy
\s@avemove \x@pix \y@pix}
\def\rmove (#1 #2){\getpos (#1 #2)\a@rgx\a@rgy
\r@elpospix \a@rgx \a@rgy
\s@avemove \x@pix \y@pix}
\def\lcir r:#1 {\coordtopix {#1}\t@pixa
\writeps {\the\t@pixa\space cr}%
\r@elupd \t@pixa \t@pixa
\r@elupd {-\t@pixa}{-\t@pixa}}
\def\fcir f:#1 r:#2 {\coordtopix {#2}\t@pixa
\writeps {\the\t@pixa\space #1 fc}%
\r@elupd \t@pixa \t@pixa
\r@elupd {-\t@pixa}{-\t@pixa}}
\def\lellip rx:#1 ry:#2 {\coordtopix {#1}\t@pixa
\coordtopix {#2}\t@pixb
\writeps {\the\t@pixa\space \the\t@pixb\space el}%
\r@elupd \t@pixa \t@pixb
\r@elupd {-\t@pixa}{-\t@pixb}}
\def\fellip f:#1 rx:#2 ry:#3 {\coordtopix {#2}\t@pixa
\coordtopix {#3}\t@pixb
\writeps {\the\t@pixa\space \the\t@pixb\space #1 fe}%
\r@elupd \t@pixa \t@pixb
\r@elupd {-\t@pixa}{-\t@pixb}}
\def\larc r:#1 sd:#2 ed:#3 {\coordtopix {#1}\t@pixa
\writeps {\the\t@pixa\space #2 #3 ar}}
\def\ifill f:#1 {\writeps {#1 fl}}
\def\lfill f:#1 {\writeps {#1 fp}}
\def\htext #1{\def\testit {#1}%
\ifx \testit\l@paren
\let\next=\h@move
\else
\let\next=\h@text
\fi
\next {#1}}
\def\rtext td:#1 #2{\def\testit {#2}%
\ifx \testit\l@paren
\let\next=\r@move
\else
\let\next=\r@text
\fi
\next td:#1 {#2}}
\def\vtext {\rtext td:90 }
\def\textref h:#1 v:#2 {\ifx #1R%
\edef\l@stuff {\hss}\edef\r@stuff {}%
\else
\ifx #1C%
\edef\l@stuff {\hss}\edef\r@stuff {\hss}%
\else
\edef\l@stuff {}\edef\r@stuff {\hss}%
\fi
\fi
\ifx #2T%
\edef\t@stuff {}\edef\b@stuff {\vss}%
\else
\ifx #2C%
\edef\t@stuff {\vss}\edef\b@stuff {\vss}%
\else
\edef\t@stuff {\vss}\edef\b@stuff {}%
\fi
\fi}
\def\avec (#1 #2){\getpos (#1 #2)\a@rgx\a@rgy
\s@etpospix \a@rgx \a@rgy
\writeps {\the\x@pix\space \the\y@pix\space (\a@type)
\the\a@lenpix\space \the\a@widpix\space av}}
\def\ravec (#1 #2){\getpos (#1 #2)\a@rgx\a@rgy
\r@elpospix \a@rgx \a@rgy
\writeps {\the\x@pix\space \the\y@pix\space (\a@type)
\the\a@lenpix\space \the\a@widpix\space av}}
\def\arrowheadsize l:#1 w:#2 {\coordtopix{#1}\a@lenpix
\coordtopix{#2}\a@widpix}
\def\arrowheadtype t:#1 {\edef\a@type{#1}}
\def\clvec (#1 #2)(#3 #4)(#5 #6)%
{\getpos (#1 #2)\a@rgx\a@rgy
\coordtopix \a@rgx\t@pixa
\advance \t@pixa by \x@segoffpix
\coordtopix \a@rgy\t@pixb
\advance \t@pixb by \y@segoffpix
\getpos (#3 #4)\a@rgx\a@rgy
\coordtopix \a@rgx\t@pixc
\advance \t@pixc by \x@segoffpix
\coordtopix \a@rgy\t@pixd
\advance \t@pixd by \y@segoffpix
\getpos (#5 #6)\a@rgx\a@rgy
\s@etpospix \a@rgx \a@rgy
\writeps {\the\t@pixa\space \the\t@pixb\space
\the\t@pixc\space \the\t@pixd\space
\the\x@pix\space \the\y@pix\space cv}}
\def\drawbb {\bsegment
\drawdim bp
\linewd 0.24
\setunitscale {\p@sfactor}
\writeps {\the\xminpix\space \the\yminpix\space mv}%
\writeps {\the\xminpix\space \the\ymaxpix\space lv}%
\writeps {\the\xmaxpix\space \the\ymaxpix\space lv}%
\writeps {\the\xmaxpix\space \the\yminpix\space lv}%
\writeps {\the\xminpix\space \the\yminpix\space lv}%
\esegment}
\def\getpos (#1 #2)#3#4{\g@etargxy #1 #2 {} \\#3#4%
\c@heckast #3%
\ifa@st
\g@etsympix #3\t@pixa
\advance \t@pixa by -\x@segoffpix
\pixtocoord \t@pixa #3%
\fi
\c@heckast #4%
\ifa@st
\g@etsympix #4\t@pixa
\advance \t@pixa by -\y@segoffpix
\pixtocoord \t@pixa #4%
\fi}
\def\getsympos (#1 #2)#3#4{\g@etargxy #1 #2 {} \\#3#4%
\c@heckast #3%
\ifa@st \else
\t@xderror {TeXdraw: invalid symbolic coordinate}%
\fi
\c@heckast #4%
\ifa@st \else
\t@xderror {TeXdraw: invalid symbolic coordinate}%
\fi}
\def\listtopix (#1)#2{\def #2{}%
\edef\l@ist {#1 }%
\m@oretrue
\loop
\expandafter\g@etitem \l@ist \\\a@rgx\l@ist
\a@pppix \a@rgx #2%
\ifx \l@ist\empty
\m@orefalse
\fi
\ifm@ore
\repeat}
\def\realmult #1#2#3{\dimen0=#1pt
\dimen2=#2\dimen0
\edef #3{\expandafter\c@lean\the\dimen2}}
\def\intdiv #1#2#3{\t@counta=#1
\t@countb=#2
\ifnum \t@countb<0
\t@counta=-\t@counta
\t@countb=-\t@countb
\fi
\t@countd=1
\ifnum \t@counta<0
\t@counta=-\t@counta
\t@countd=-1
\fi
\t@countc=\t@counta  \divide \t@countc by \t@countb
\t@counte=\t@countc  \multiply \t@counte by \t@countb
\advance \t@counta by -\t@counte
\t@counte=-1
\loop
\advance \t@counte by 1
\ifnum \t@counte<16
\multiply \t@countc by 2
\multiply \t@counta by 2
\ifnum \t@counta<\t@countb \else
\advance \t@countc by 1
\advance \t@counta by -\t@countb
\fi
\repeat
\divide \t@countb by 2
\ifnum \t@counta<\t@countb
\advance \t@countc by 1
\fi
\ifnum \t@countd<0
\t@countc=-\t@countc
\fi
\dimen0=\t@countc sp
\edef #3{\expandafter\c@lean\the\dimen0}}
\def\coordtopix #1#2{\dimen0=#1\d@dim
\dimen2=\d@sc\dimen0
\t@counta=\dimen2
\t@countb=\s@ppix
\divide \t@countb by 2
\ifnum \t@counta<0
\advance \t@counta by -\t@countb
\else
\advance \t@counta by \t@countb
\fi
\divide \t@counta by \s@ppix
#2=\t@counta}
\def\pixtocoord #1#2{\t@counta=#1%
\multiply \t@counta by \s@ppix
\dimen0=\d@sc\d@dim
\t@countb=\dimen0
\intdiv \t@counta \t@countb #2}
\def\pixtodim #1#2{\t@countb=#1%
\multiply \t@countb by \s@ppix
#2=\t@countb sp\relax}
\def\pixtobp #1#2{\dimen0=\p@sfactor pt
\t@counta=\dimen0
\multiply \t@counta by #1%
\ifnum \t@counta < 0
\advance \t@counta by -32768
\else
\advance \t@counta by 32768
\fi
\divide \t@counta by 65536
\edef #2{\the\t@counta}}
\newcount\t@counta    \newcount\t@countb
\newcount\t@countc    \newcount\t@countd
\newcount\t@counte
\newcount\t@pixa      \newcount\t@pixb
\newcount\t@pixc      \newcount\t@pixd
\newdimen\t@xpos      \newdimen\t@ypos
\newcount\xminpix      \newcount\xmaxpix
\newcount\yminpix      \newcount\ymaxpix
\newcount\a@lenpix     \newcount\a@widpix
\newcount\x@pix        \newcount\y@pix
\newcount\x@segoffpix  \newcount\y@segoffpix
\newcount\x@savepix    \newcount\y@savepix
\newcount\s@ppix
\newcount\d@bs
\newcount\t@xdnum
\global\t@xdnum=0
\newbox\t@xdbox
\newwrite\drawfile
\newif\ifm@pending
\newif\ifp@ath
\newif\ifa@st
\newif\ifm@ore
\newif \ift@extonly
\newif\ifp@osinit
\newtoks\everytexdraw
\def\l@paren{(}
\def\a@st{*}
\catcode`\%=12
\def\p@b {
\catcode`\%=14
\catcode`\{=12  \catcode`\}=12  \catcode`\u=1 \catcode`\v=2
\def\l@br u{v  \def\r@br u}v
\catcode `\{=1  \catcode`\}=2   \catcode`\u=11 \catcode`\v=11
{\catcode`\p=12 \catcode`\t=12
\gdef\c@lean #1pt{#1}}
\def\sppix#1/#2 {\dimen0=1#2 \s@ppix=\dimen0
\t@counta=#1%
\divide \t@counta by 2
\advance \s@ppix by \t@counta
\divide \s@ppix by #1%
\t@counta=\s@ppix
\multiply \t@counta by 65536
\advance \t@counta by 32891
\divide \t@counta by 65782
\dimen0=\t@counta sp
\edef\p@sfactor {\expandafter\c@lean\the\dimen0}}
\def\g@etargxy #1 #2 #3 #4\\#5#6{\def #5{#1}%
\ifx #5\empty
\g@etargxy #2 #3 #4 \\#5#6
\else
\def #6{#2}%
\def\next {#3}%
\ifx \next\empty \else
\t@xderror {TeXdraw: invalid coordinate}%
\fi
\fi}
\def\c@heckast #1{\expandafter
\c@heckastll #1\\}
\def\c@heckastll #1#2\\{\def\testit {#1}%
\ifx \testit\a@st
\a@sttrue
\else
\a@stfalse
\fi}
\def\g@etsympix #1#2{\expandafter
\ifx \csname #1\endcsname \relax
\t@xderror {TeXdraw: undefined symbolic coordinate}%
\fi
#2=\csname #1\endcsname}
\def\s@etcsn #1#2{\expandafter
\xdef\csname#1\endcsname {#2}}
\def\g@etitem #1 #2\\#3#4{\edef #4{#2}\edef #3{#1}}
\def\a@pppix #1#2{\edef\next {#1}%
\ifx \next\empty \else
\coordtopix {#1}\t@pixa
\ifx #2\empty
\edef #2{\the\t@pixa}%
\else
\edef #2{#2 \the\t@pixa}%
\fi
\fi}
\def\s@etpospix #1#2{\coordtopix {#1}\x@pix
\advance \x@pix by \x@segoffpix
\coordtopix {#2}\y@pix
\advance \y@pix by \y@segoffpix
\u@pdateminmax \x@pix \y@pix}
\def\r@elpospix #1#2{\coordtopix {#1}\t@pixa
\advance \x@pix by \t@pixa
\coordtopix {#2}\t@pixa
\advance \y@pix by \t@pixa
\u@pdateminmax \x@pix \y@pix}
\def\r@elupd #1#2{\t@counta=\x@pix
\advance\t@counta by #1%
\t@countb=\y@pix
\advance\t@countb by #2%
\u@pdateminmax \t@counta \t@countb}
\def\u@pdateminmax #1#2{\ifnum #1>\xmaxpix
\global\xmaxpix=#1%
\fi
\ifnum #1<\xminpix
\global\xminpix=#1%
\fi
\ifnum #2>\ymaxpix
\global\ymaxpix=#2%
\fi
\ifnum #2<\yminpix
\global\yminpix=#2%
\fi}
\def\s@avemove #1#2{\x@savepix=#1\y@savepix=#2%
\m@pendingtrue
\ifp@osinit \else
\global\p@osinittrue
\global\xminpix=\x@savepix \global\yminpix=\y@savepix
\global\xmaxpix=\x@savepix \global\ymaxpix=\y@savepix
\fi}
\def\f@lushmove {\global\p@osinittrue
\ifm@pending
\writetx {\the\x@savepix\space \the\y@savepix\space mv}%
\m@pendingfalse
\p@athfalse
\fi}
\def\f@lushbs {\loop
\ifnum \d@bs>0
\writetx {bs}%
\global\advance \d@bs by -1
\repeat}
\def\h@move #1#2 #3)#4{\move (#2 #3)%
\h@text {#4}}
\def\h@text #1{\pixtodim \x@pix \t@xpos
\pixtodim \y@pix \t@ypos
\vbox to 0pt{\normalbaselines
\t@stuff
\kern -\t@ypos
\hbox to 0pt{\l@stuff
\kern \t@xpos
\hbox {#1}%
\kern -\t@xpos
\r@stuff}%
\kern \t@ypos
\b@stuff\relax}}
\def\r@move td:#1 #2#3 #4)#5{\move (#3 #4)%
\r@text td:#1 {#5}}
\def\r@text td:#1 #2{\vbox to 0pt{\pixtodim \x@pix \t@xpos
\pixtodim \y@pix \t@ypos
\kern -\t@ypos
\hbox to 0pt{\kern \t@xpos
\rottxt {#1}{\z@sb {#2}}%
\hss}%
\vss}}
\def\z@sb #1{\vbox to 0pt{\normalbaselines
\t@stuff
\hbox to 0pt{\l@stuff \hbox {#1}\r@stuff}%
\b@stuff}}
\ifx \rotatebox\@undefined
\def\rottxt #1#2{\bgroup
#2%
\egroup}
\else
\let\rottxt=\rotatebox
\fi
\ifx \t@xderror\@undefined
\let\t@xderror=\errmessage
\fi
\def\t@exdrawdef {\sppix 300/in
\drawdim in
\edef\u@nitsc {1}%
\setsegscale 1
\arrowheadsize l:0.16 w:0.08
\arrowheadtype t:T
\textref h:L v:B }
\ifx \includegraphics\@undefined
\def\t@xdinclude [#1,#2][#3,#4]#5{%
\begingroup
\message {<#5>}%
\leavevmode
\t@counta=-#1%
\t@countb=-#2%
\setbox0=\hbox{%
\includegraphics{#5}}%
\t@ypos=#4 bp%
\advance \t@ypos by -#2 bp%
\t@xpos=#3 bp%
\advance \t@xpos by -#1 bp%
\dp0=0pt \ht0=\t@ypos  \wd0=\t@xpos
\box0%
\endgroup}
\else
\let\t@xdinclude=\includegraphics
\fi
\def\writeps #1{\f@lushbs
\f@lushmove
\p@athtrue
\writetx {#1}}
\def\writetx #1{\ift@extonly
\global\t@extonlyfalse
\t@xdpsfn \p@sfile
\t@dropen \p@sfile
\fi
\w@rps {#1}}
\def\w@rps #1{\immediate\write\drawfile {#1}}
\def\t@xdpsfn #1{%
\global\advance \t@xdnum by 1
\ifnum \t@xdnum<10
\xdef #1{\jobname.ps\the\t@xdnum}
\else
\xdef #1{\jobname.p\the\t@xdnum}
\fi
}
\def\t@dropen #1{%
\immediate\openout\drawfile=#1%
\w@rps {\p@b PS-Adobe-3.0 EPSF-3.0}%
\w@rps {\p@p BoundingBox: (atend)}%
\w@rps {\p@p Title: TeXdraw drawing: #1}%
\w@rps {\p@p Pages: 1}%
\w@rps {\p@p Creator: \TeXdrawId}%
\w@rps {\p@p CreationDate: \the\year/\the\month/\the\day}%
\w@rps {50 dict begin}%
\w@rps {/mv {stroke moveto} def}%
\w@rps {/lv {lineto} def}%
\w@rps {/st {currentpoint stroke moveto} def}%
\w@rps {/sl {st setlinewidth} def}%
\w@rps {/sd {st 0 setdash} def}%
\w@rps {/sg {st setgray} def}%
\w@rps {/bs {gsave} def /es {stroke grestore} def}%
\w@rps {/fl \l@br gsave setgray fill grestore}%
\w@rps    { currentpoint newpath moveto\r@br\space def}%
\w@rps {/fp {gsave setgray fill grestore st} def}%
\w@rps {/cv {curveto} def}%
\w@rps {/cr \l@br gsave currentpoint newpath 3 -1 roll 0 360 arc}%
\w@rps    { stroke grestore\r@br\space def}%
\w@rps {/fc \l@br gsave setgray currentpoint newpath}%
\w@rps    { 3 -1 roll 0 360 arc fill grestore\r@br\space def}%
\w@rps {/ar {gsave currentpoint newpath 5 2 roll arc stroke grestore} def}%
\w@rps {/el \l@br gsave /svm matrix currentmatrix def}%
\w@rps    { currentpoint translate scale newpath 0 0 1 0 360 arc}%
\w@rps    { svm setmatrix stroke grestore\r@br\space def}%
\w@rps {/fe \l@br gsave setgray currentpoint translate scale newpath}%
\w@rps    { 0 0 1 0 360 arc fill grestore\r@br\space def}%
\w@rps {/av \l@br /hhwid exch 2 div def /hlen exch def}%
\w@rps    { /ah exch def /tipy exch def /tipx exch def}%
\w@rps    { currentpoint /taily exch def /tailx exch def}%
\w@rps    { /dx tipx tailx sub def /dy tipy taily sub def}%
\w@rps    { /alen dx dx mul dy dy mul add sqrt def}%
\w@rps    { /blen alen hlen sub def}%
\w@rps    { gsave tailx taily translate dy dx atan rotate}%
\w@rps    { (V) ah ne {blen 0 gt {blen 0 lineto} if} {alen 0 lineto} ifelse}%
\w@rps    { stroke blen hhwid neg moveto alen 0 lineto blen hhwid lineto}%
\w@rps    { (T) ah eq {closepath} if}%
\w@rps    { (W) ah eq {gsave 1 setgray fill grestore closepath} if}%
\w@rps    { (F) ah eq {fill} {stroke} ifelse}%
\w@rps    { grestore tipx tipy moveto\r@br\space def}%
\w@rps {\p@sfactor\space \p@sfactor\space scale}%
\w@rps {1 setlinecap 1 setlinejoin}%
\w@rps {3 setlinewidth [] 0 setdash}%
\w@rps {0 0 moveto}%
}
\def\t@drclose {%
\bgroup
\w@rps {stroke end showpage}%
\w@rps {\p@p Trailer:}%
\pixtobp \xminpix \l@lxbp  \pixtobp \yminpix \l@lybp
\pixtobp \xmaxpix \u@rxbp  \pixtobp \ymaxpix \u@rybp
\w@rps {\p@p BoundingBox: \l@lxbp\space \l@lybp\space
\u@rxbp\space \u@rybp}%
\w@rps {\p@p EOF}%
\egroup
\immediate\closeout\drawfile
}
\catcode`\@=\catamp

\def\ldreieck{\bsegment
  \rlvec(0.866025403784439 .5) \rlvec(0 -1)
  \rlvec(-0.866025403784439 .5)  
  \savepos(0.866025403784439 -.5)(*ex *ey)
        \esegment
  \move(*ex *ey)
        }
\def\rdreieck{\bsegment
  \rlvec(0.866025403784439 -.5) \rlvec(-0.866025403784439 -.5)  \rlvec(0 1)
  \savepos(0 -1)(*ex *ey)
        \esegment
  \move(*ex *ey)
        }
\def\rhombus{\bsegment
  \rlvec(0.866025403784439 .5) \rlvec(0.866025403784439 -.5) 
  \rlvec(-0.866025403784439 -.5)  \rlvec(0 1)        
  \rmove(0 -1)  \rlvec(-0.866025403784439 .5) 
  \savepos(0.866025403784439 -.5)(*ex *ey)
        \esegment
  \move(*ex *ey)
        }
\def\RhombusA{\bsegment
  \rlvec(0.866025403784439 .5) \rlvec(0.866025403784439 -.5) 
  \rlvec(-0.866025403784439 -.5) \rlvec(-0.866025403784439 .5) 
  \savepos(0.866025403784439 -.5)(*ex *ey)
        \esegment
  \move(*ex *ey)
        }
\def\RhombusB{\bsegment
  \rlvec(0.866025403784439 .5) \rlvec(0 -1)
  \rlvec(-0.866025403784439 -.5) \rlvec(0 1) 
  \savepos(0 -1)(*ex *ey)
        \esegment
  \move(*ex *ey)
        }
\def\RhombusC{\bsegment
  \rlvec(0.866025403784439 -.5) \rlvec(0 -1)
  \rlvec(-0.866025403784439 .5) \rlvec(0 1) 
  \savepos(0.866025403784439 -.5)(*ex *ey)
        \esegment
  \move(*ex *ey)
        }
\def\hdSchritt{\bsegment
  \lpatt(.05 .13)
  \rlvec(0.866025403784439 -.5) 
  \savepos(0.866025403784439 -.5)(*ex *ey)
        \esegment
  \move(*ex *ey)
        }
\def\vdSchritt{\bsegment
  \lpatt(.05 .13)
  \rlvec(0 -1) 
  \savepos(0 -1)(*ex *ey)
        \esegment
  \move(*ex *ey)
        }

\def\ringerl(#1 #2){\move(#1 #2)\fcir f:0 r:.07}

\def\hex{\bsegment
	\rlvec(1 0)  \rlvec(.5 -.866025403784439) \rlvec(-.5 -.866025403784439)
	\rlvec(-1 0) \rlvec(-.5 .866025403784439) \rlvec(.5 .866025403784439)
	\savepos(1.5 -.866025403784439)(*ex *ey)
	 \esegment
	\move(*ex *ey)
}

\def\shex{\bsegment 
	\lpatt(.05 .13)
	\rlvec(1 0)  \rlvec(.5 -.866025403784439) \rlvec(-.5 -.866025403784439)
	\rlvec(-1 0) \rlvec(-.5 .866025403784439) \rlvec(.5 .866025403784439)
	\savepos(1.5 -.866025403784439)(*ex *ey)
	 \esegment
	\move(*ex *ey)
}

\def\({\left(}
\def\){\right)}
\def\[{\left[}
\def\]{\right]}
\def\fl#1{\left\lfloor#1\right\rfloor}
\def\cl#1{\left\lceil#1\right\rceil}
\def\odd{\widetilde}

\begin{section}{Introduction} \label{introsec}
The interest in rhombus tilings has emerged from the enumeration of
plane partitions in a given box.
The connection comes from looking at
the stacks of cubes of a plane partition 
from the right angle and projecting the
picture to the plane. Then the box becomes a hexagon, where opposite
sides are equal, and the cubes become a rhombus tiling of the hexagon 
where the rhombi consist of two equilateral triangles (cf\@.
\cite{David-Tomei}).
The number of plane partitions in a given box was first computed by
MacMahon        
\cite[Sec.~429, $q\to1$, proof in Sec.~494]{MM}. Therefore:

\smallskip
{\it The number of all rhombus tilings of a hexagon
with sides $a,b,c,a,b,c$ equals}
\begin{equation} \label{box}
B(a,b,c)=\prod _{i=1} ^{a}\prod _{j=1} ^{b}\prod _{k=1} ^{c}\frac {i+j+k-1}
{i+j+k-2}.
\end{equation}
(The form of the expression is due to Macdonald.)

In \cite{Propp}, Propp proposed several problems regarding
``incomplete'' hexagons, i.e., hexagons, where certain triangles are
missing. In particular, Problem~4 of \cite{Propp} asks for a formula     
for the number of 
rhombus tilings of a regular hexagon, where two of the six central
triangles are missing. We treat the case of the two triangles
lying on the symmetry axis and touching in one vertex (see
Figure~\ref{fliegefi}). The other case has been solved in  
\cite{FuKrAD}. We prove the following two theorems. 

\begin{theorem} \label{th1}
The number of rhombus tilings of a hexagon with sides
$n, n, 2m, n, n, 2m$ and two missing triangles on the horizontal 
symmetry axis sharing the $(s+1)$--th vertex on the axis (see
Figure~\ref{fliegefi}) equals
$$
\frac { (2m-1) \binom{2m-2}{m-1} \binom{2n-2s}{n-s} \binom{2s}{s}} 
{\binom{2m+2n}{m+n}}
\prod _{i=1} ^{n}\prod _{j=1} ^{n}
\prod _{k=1} ^{2m}{\frac {i+j+k-1} {i+j+k-2}}.
$$
\end{theorem}
\begin{theorem}
\label{th2}
The number of rhombus tilings of a hexagon with sides
$n,n,2m+1,$ $n,n,2m+1$ and two missing triangles on the symmetry axis
sharing the $s$--th vertex on the axis equals
$$
\frac {(2m+1) \binom{2m}{m} \binom{2n-2s}{n-s} \binom{2s-2}{s-1}}
{\binom{2m+2n}{m+n}} 
\prod _{i=1} ^{n} \prod _{j=1} ^{n} 
\prod _{k=1} ^{2m+1}{\frac {i+j+k-1} {i+j+k-2}}.
$$
\end{theorem}

The following corollary is easily derived using Stirling's
approximation formula.

\begin{cor}
The proportion of rhombus tilings of a hexagon with sides
$\alpha t, \alpha t, \beta t,\alpha t,\alpha t, \beta t$ and 
two missing triangles on the horizontal symmetry axis touching
the $(\gamma t)$-th vertex on the axis in the number 
of all rhombus tilings of the hexagon with sidelengths
$\alpha t, \alpha t, \beta t,\alpha t,\alpha t, \beta t$ 
(given by 
$B(\alpha t, \alpha t, \beta t,\alpha t,\alpha t, \beta t)$ in
\eqref{box})
is asymptotically equal to
$$\frac {1} {4 \pi} \sqrt {\frac {\beta (2\alpha + \beta)} {\gamma
(\alpha - \gamma)}}.$$
\end{cor}

This expression can attain arbitrary large values if $\alpha$ is
close to $\gamma$. It equals $\frac
{\sqrt 3} {2\pi}$ (which is approximately 0.28)
for $\alpha =\beta = 2 \gamma$, which corresponds to the case of a
regular hexagon with two missing triangles touching the center.
In comparison, in the other case of Problem~4 of \cite{Propp},
the case of a fixed rhombus on the symmetry axis, the analogous 
proportion must clearly be smaller than 1 and equals approximately 
$\frac {1} {3}$ if the central rhombus is missing (see \cite{FuKrAD}).
  
The rest of the paper is devoted to the proof of Theorems~\ref{th1}
and \ref{th2}. The
main ingredients are the matchings factorization theorem by 
M.~Ciucu \cite{Ci1}, nonintersecting lattice paths, and two
determinant evaluations, the latter constituting the most difficult
part of the proof. An outline of the proof is given in the next
section. The details are filled in in the subsequent sections.
\vskip3cm
\begin{figure}
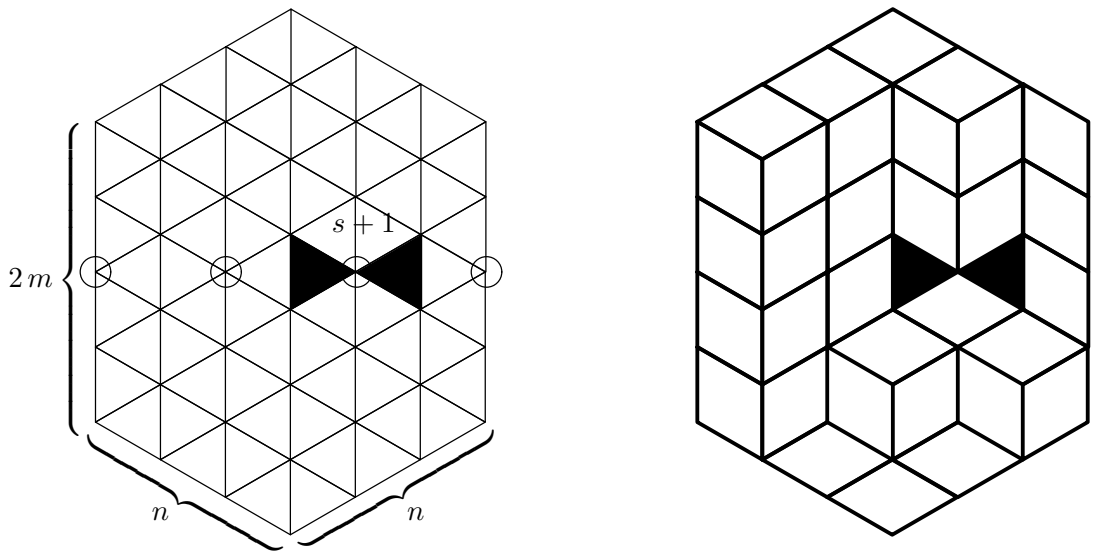

\centertexdraw{
\drawdim truecm \linewd.02
\rhombus \rhombus \rhombus \ldreieck
\move(-.866025 -.5) 
\rhombus \rhombus \rhombus \rhombus \ldreieck
\move(-1.73205 -1) 
\rhombus \rhombus \rhombus \rhombus \rhombus \ldreieck
\move(-1.73205 -1)
\rdreieck \rhombus \rhombus \rhombus \rhombus \rhombus \ldreieck
\move(-1.73205 -2)
\rdreieck \rhombus \rhombus \rhombus \rhombus \rhombus
\move(-1.73205 -3)
\rdreieck \rhombus \rhombus \rhombus \rhombus
\move(-1.73205 -4)
\rdreieck \rhombus \rhombus \rhombus

\move(.866025 -2.5)
\bsegment
 \rlvec(0.866025403784439 -.5) \rlvec(-0.866025403784439 -.5) 
 \lfill f:0
\esegment
\move(1.73205 -3)
\bsegment
 \rlvec(0.866025403784439 -.5) \rlvec(0 1) 
 \lfill f:0
\esegment

\move(0 -3)
\lcir r:.2

\move(-1.73205 -3)
\lcir r:.2

\move(1.73205 -3)
\lcir r:.2

\move(3.464 -3)
\lcir r:.2

\htext(1.4 -2.5){$s+1$}
\rtext td:0 (-2.9 -3.2) {$2\, m \left\{ \vbox{\vskip2.2cm} \right.$ }
\rtext td:60 (-.6 -6.2) {$\left\{ \vbox{\vskip1.6cm}\right.$ }
\rtext td:-60 (2.1 -5.8) {$\left.\vbox{\vskip1.6cm}\right\}$}
\htext(-1 -6.3) {$n$}
\htext(2.4 -6.3) {$n$}

\move(8 0)
\bsegment
\linewd.01

\move(.866025 -2.5)
\bsegment
 \rlvec(0.866025403784439 -.5) \rlvec(-0.866025403784439 -.5) 
 \lfill f:0
  \esegment
\move(1.73205 -3)
\bsegment
 \rlvec(0.866025403784439 -.5) \rlvec(0 1) 
 \lfill f:0
  \esegment

\linewd.05
\move(0 0)
\RhombusA \RhombusA \RhombusB \RhombusB
\move(-.866025 -.5)
\RhombusA \RhombusB \RhombusB \RhombusB \RhombusA \RhombusB \RhombusA
\move(-1.73205 -1)
\RhombusA \RhombusB \RhombusB \RhombusB \RhombusB \RhombusA \RhombusA
\move(-1.73205 -1)
\RhombusC 
\move(-1.73205 -2)
\RhombusC 
\move(-1.73205 -3)
\RhombusC 
\move(-1.73205 -4)
\RhombusC
\move(0.866025 -.5) \RhombusC \move(.866025 -1.5) \RhombusC
\move(2.598075 -.5) \RhombusC
\move(2.598075 -1.5) \RhombusC
\move(2.598075 -2.5) \RhombusC
\move(.866025 -3.5) \RhombusA \RhombusA \RhombusB
\esegment
}
\caption{ \label{fliegefi} A hexagon with sides $n,n,2m,n,n,2m$ and 
missing triangles in position $s+1$, where $m=2,n=3,s=2$, and a 
rhombus tiling.}
\end{figure}
\end{section}
\begin{section}{Outline of the proofs of Theorems~\ref{th1} and
\ref{th2}} \label{outlsec}
{\bf Outline of the proof of Theorem~\ref{th1}:}

\smallskip
{\it Step 1: It suffices to compute the number of perfect matchings
of two hexagonal graphs $G^+$ and $G^-$ (see Section~\ref{cutsec}).}

We use the fact that every rhombus tiling corresponds to a
perfect matching of the (inner) dual graph (see Figure~\ref{hexfi}a).
This graph has reflective symmetry, so the matchings 
factorization theorem by M.~Ciucu \cite{Ci1}
is applicable (see Lemma~\ref{fact}). This theorem expresses
the number of perfect matchings of a graph as product of $2^{n-1}$
and the numbers of perfect matchings of two smaller graphs $G^+$ and
$G^-$ (see Lemma~\ref{split} and Figure~\ref{hexfi}), which are 
roughly the two halves of 
the original graph. 
The remaining task is to count the numbers of perfect matchings of $G^+$ and
$G^-$.

\smallskip
{\it Step 2: The numbers of perfect matchings of $G^+$ and $G^-$
equal the numbers of rhombus tilings of two regions $R^+$ and $R^-$
consisting of triangles, respectively 
(see Sections~\ref{upsec} and \ref{lowsec}).}
We convert the graphs $G^+$ and $G^-$ back to regions $R^+$ and $R^-$ 
consisting of triangles (see Figure~\ref{uplowfi}), so that we have to count 
rhombus tilings again.

\smallskip
{\it Step 3: The numbers of rhombus tilings of $R^+$ and $R^-$ are 
certain determinants (see Sections~\ref{upsec} and \ref{lowsec}).}

The rhombus tilings are in bijection with certain families of 
nonintersecting lattice paths (see Figures~\ref{upfi} and \ref{lowfi}).
Application of the main result of nonintersecting lattice paths
expresses the desired numbers as determinants (see Lemmas~\ref{plusdet}
and \ref{minusdet}).

\smallskip
{\it Step 4: Evaluation of the determinant corresponding to $G^+$ (see
Section~\ref{upsec}).}

The determinant corresponding to $G^+$ is evaluated using a
lemma by Krattenthaler (see Lemmas~\ref{det} and \ref{g+}).

\smallskip
{\it Step 5: Evaluation of the determinant corresponding to $G^-$.}

We pull factors out of the determinant, 
so that we obtain a determinant whose entries are polynomials in $m$
(see the proof of Lemma~\ref{minusdet}).
Then, in Lemma~\ref{g-}, this determinant is evaluated by using the
``identification of factors'' method, as explained in
\cite[Sec.~2]{KratBI}. The corresponding details are the subject of
Sections~\ref{halbsec} -- \ref{degsec}. 

\smallskip
{\it Step 6: A combination of the previous steps proves Theorem~\ref{th1}.}

We substitute the results of Lemmas~\ref{g+}, \ref{minusdet} and 
\ref{g-} in Lemma~\ref{fact} and obtain the following expression 
for the number of rhombus tilings of our original hexagon,
\begin{multline*}
2^{\binom{n}{2}-1}\frac {\left(\h(n)\right)^2(2n-2s-1)!!(2s-1)!!}
{\h(2n)(n-s)!s!}\\
\times\prod _{2\le i\le j\le n} ^{}{(2m+2j-i)}
\prod _{k=1} ^{n-2}{\left( m+k+\frac{1}{2}\right)^{\min(k,n-1-k)} 
\prod_{k=0} ^{n}{(m+k)^{\min(k+1,n-k+1)}}}.
\end{multline*}
Here, $\h(n)$ stands for $\prod_{i=0}^{n-1}{i!}$.
This can easily be transformed to the expression in
Theorem~\ref{th1}, so the
proof of Theorem~\ref{th1} is complete. 

\smallskip
{\bf Outline of the proof of Theorem~\ref{th2} (see
Section~\ref{oddsec}):}

\smallskip
{\it Step $1'$: It is enough to count the rhombus tilings of the two
regions $\odd{R^+}$ and $\odd{R^-}$ shown in Figure~\ref{oddfi}.}

By the analogue of Step~1, we have to count the numbers of perfect
matchings of two graphs $\odd{G^+}$ and $\odd{G^-}$. 
These numbers equal the numbers of rhombus 
tilings of the regions $\odd{R^+}$ and $\odd{R^-}$ of triangles shown in 
Figure~\ref{oddfi} (compare Step~2 in the previous proof).

\smallskip
{\it Step $2'$: The number of rhombus tilings of $\odd{R^+}$.}

The number of rhombus tilings of $\odd{R^+}(n,m)$ equals the number
of rhombus tilings of $R^+(n-1,m)$, so it can be found using
Lemma~\ref{g+} (see equation~\eqref{odd+}).

\smallskip
{\it Step $3'$: The number of rhombus tilings of $\odd{R^-}$.}

The number of rhombus tilings of $\odd{R^-}(n,m,s)$ 
can be expressed as a determinant analogously to Step 2 in the
previous proof.
The determinant equals the determinant counting
the number of rhombus tilings of $R^-(n-1,m+1,s-1)$, so 
Lemma~\ref{g-} can be used for the evaluation (see equation~\eqref{odd-}). 

\smallskip
{\it Step $4'$: A combination of the previous steps proves Theorem~\ref{th2}.}
We substitute the results of
equations~\eqref{odd+} and \eqref{odd-} in \eqref{facto}. We get 
\begin{align*} 
M(\odd{G})&=2^{n-1}M(\odd{R^+})M(\odd{R^-})\\
&=\frac {2^{\binom{n-1}{2}-1}\h(n-1)\h(n+1)(2n-2s-1)!!(2s-3)!!}
{(n-s)! (s-1)! \prod _{j=0} ^{n}{(2j)!}\prod _{i=0}
^{n-2}{(2i+1)!}}\nonumber\\
&\times
\prod _{k=2} ^{n-2}{\left( m+k+\frac{1}{2}\right)^{\min(k-1,n-k-1)} 
\prod_{k=1} ^{n}{(m+k)^{\min(k,n-k+1)}}}
\prod _{2\le i \le j \le n+1} ^{}{(2m+2j-i)},\nonumber
\end{align*}
which can easily be transformed to the expression in Theorem~\ref{th2}.
\end{section}
\begin{section}{Breaking the hexagon in two parts} \label{cutsec}
We start the proof of Theorem~\ref{th1} by forming the inner dual
of the given hexagon.
I\@.e\@., we replace every triangle by a vertex and connect vertices
corresponding to adjacent triangles. Thus, we get a hexagonal graph, whose
perfect matchings correspond to rhombus tilings of the original
hexagon (see Figure~\ref{hexfi}a). 

Now we use a theorem by M.~Ciucu (see \cite{Ci1}) to, roughly
speaking, break the
hexagonal graph into two halves. This theorem is described as follows.
Let $G$ be a graph with reflective symmetry, which splits into two 
parts after removal of the vertices of the symmetry axis.
Label the vertices on the symmetry axis
$a_1,b_1,a_2,b_2,\dots,a_{2l},b_{2l}$ from left to right.
If $G$ is bipartite, we can colour
the vertices of the graph black and white subject to the conditions
that $a_1$ is white and no two adjacent vertices are of the same
colour.

Then we delete all edges connecting white $a$--vertices and black
$b$--vertices to the upper half and all edges connecting black
$a$--vertices and white $b$--vertices to the lower half.
If we divide by two 
all weights of edges lying on the symmetry axis, the graph splits
into two parts $G^+$ and $G^-$. Now we can state the matchings
factorization theorem from~\cite{Ci1}.
\begin{lemma} \label{fact}
Let $G$ be a planar bipartite weighted, symmetric 
graph, which splits into two parts after removal of the
vertices of the symmetry axis.
Then $$M(G)=2^{l(G)}M(G^+)M(G^-),$$ where
$M(H)$ denotes the weighted count
of perfect matchings of the graph $H$ and $G^\pm$ denote the upper
and lower half of $G$ as described above.
$2\, l(G)$ is the number of vertices on the symmetry axis.
\end{lemma}
We apply Lemma~\ref{fact} to our hexagonal graph, exemplified in
Figure~\ref{hexfi}a, with respect to the horizontal symmetry axis.
In our case $l(G)=n-1$, since two vertices correspond to the removed
triangles. $G^+$ and $G^-$ are shown in Figure~\ref{hexfi}b. 
Thus we get the following lemma.
\begin{lemma} \label{split}
The number of rhombus tilings of a hexagon with sides
$n, n, 2m, n, n, 2m$ and two missing triangles on the symmetry axis
sharing the $(s+1)$--th vertex on the axis equals
$$2^{n-1} M(G^+) M(G^-),$$
where $G^+$ and $G^-$ are formed by the above procedure, as
exemplified in Figure~\ref{hexfi}b. 
\end{lemma}

$M(G^+)$ and $M(G^-)$ are computed in the following sections.
\end{section}
\begin{section}{The matchings count for the upper half} \label{upsec}
In this section we evaluate $M(G^+)$.
We start by expressing $M(G^+)$ as the following determinant.

\begin{lemma}\label{plusdet}
$$M(G^+)=\det_{1\le i,j\le n}\(\binom{m+j-1}{m-j+i}\).$$
\end{lemma}
\begin{proof}
First, we convert $G^+$ back to the corresponding region of
triangles, $R^+$ say 
(see Figure~\ref{uplowfi}a), so that perfect matchings of
$G^+$ correspond bijectively to the rhombus tilings of $R^+$.
Thus, we have to count rhombus tilings of $R^+$.
The next step is converting rhombus tilings to families of
nonintersecting lattice paths, i\@.e\@. lattice paths which have 
no common vertices. The reader should consult Figure~\ref{upfi}, 
while reading the following passage. Given a rhombus tiling of $R^+$,
the lattice paths start on the centers of upper left diagonal edges 
(lying on one of the sides of length $n$). They end on the lower right 
edges parallel to the starting edges. The paths are generated by 
connecting the center of the respective edge with the center of the 
edge lying opposite in the rhombus. This process is iterated using the 
new edge and the second rhombus it bounds. It terminates on the lower 
right boundary edges. It is obvious that paths starting at different 
points have no common vertices, i\@.e\@., are nonintersecting.
Furthermore, an arbitrary family of nonintersecting paths from
the set of the upper left edges to the set of the lower right edges
lies completely inside $R^+$ and can be converted back to a
tiling (see Figure~\ref{upfi}a).

Then we transform the picture to ``orthogonal'' paths with positive
horizontal and negative vertical steps of unit length (see
Figure~\ref{upfi}b,c).
Let the starting points of the paths be denoted by $P_1, P_2, \dots ,
P_n$ and the end points by $Q_1, Q_2, \dots ,Q_n$. We can easily
write down the coordinates of the starting points and the end points: 
\begin{align*}
P_i&=(i-1,i+m-1)\quad &\text {for }i=1,\dots,n,\\
Q_j&=(2j-2,j-1)\quad &\text {for }j=1,\dots,n.
\end{align*}
Next we apply the main result for nonintersecting lattice paths 
\cite[Cor.2]{gv} (see also \cite[Theorem~1.2]{StemAE}).
This theorem says that the weighted count of families of nonintersecting
lattice paths, with path $i$ running from $P_i$ to $Q_i$, is the
determinant of the matrix with $(i,j)$-entry the 
weight $\mathcal{P}(P_i \rightarrow Q_j)$ 
of lattice paths running from $P_i$ to $Q_j$, provided
that every two paths $P_i \to Q_j$ and $P_k \to Q_l$ 
have a common vertex if $i<j$ and $k>l$.
It is easily checked that our sets of starting and end points 
meet the required conditions.
\begin{figure}
\centertexdraw{
\drawdim cm
\setunitscale.5773503
\linewd.04
\shex \shex \shex
\move(-1.5 -.866025)
\shex \shex \shex
\move(-3 -1.73205)
\shex \shex 
\move(-3 -3.4641)
\shex \shex
\move(-3 -5.19615)
\shex \shex \shex \shex
\move(-3 -6.9282)
\shex \shex \shex
\move(1.5 -6.062168)
\shex \shex
\move(4 -3.4641)
\lpatt(.05 .13)
\rlvec(.5 -.866025)
\rlvec(-.5 -.866025)
\rlvec(.5 -.866025)
\rlvec(-.5 -.866025)
\lpatt(1 0) \linewd.08
\move(0 0) \rlvec(1 0)
\rmove(.5 -.866025) \rlvec(1 0)
\rmove(.5 -.866025) \rlvec(-.5 -.866025)
\rmove(.5 -.866025) \rlvec(-.5 -.866025)
\move(-1.5 -.866025) \rlvec(1 0)
\rmove(.5 -.866025) \rlvec(-.5 -.866025)
\rmove(.5 -.866025) \rlvec(-.5 -.866025)
\rmove(.5 -.866025) \rlvec(-.5 -.866025)
\rmove(.5 -.866025) \rlvec(1 0)
\rmove(.5 -.866025)\rlvec(-.5 -.866025)
\rmove(.5 -.866025) \rlvec(1 0)

\move(-3 -1.73205)
\rlvec(1 0)
\rmove(.5 -.866025) \rlvec(-.5 -.866025)
\rmove(.5 -.866025) \rlvec(-.5 -.866025)
\rmove(.5 -.866025) \rlvec(-.5 -.866025)
\rmove(.5 -.866025) \rlvec(-.5 -.866025)
\rmove(.5 -.866025) \rlvec(1 0)
\rmove(.5 -.866025) \rlvec(1 0)
\move(-3.5 -2.598075)
\rlvec(.5 -.866025) \rmove(-.5 -.866025)
\rlvec(.5 -.866025) \rmove(-.5 -.866025)
\rlvec(.5 -.866025) \rmove(-.5 -.866025)
\rlvec(.5 -.866025) \rmove(-.5 -.866025)

\move(4 -1.73205)
\rlvec(.5 -.866025) \rmove(-.5 -.866025)
\rlvec(.5 -.866025) \rmove(-.5 -.866025)
\rlvec(.5 -.866025) \rmove(-.5 -.866025)
\move(1 -1.73205)
\rlvec(.5 -.866025) \rmove(-.5 -.866025)
\rlvec(.5 -.866025) \rmove(-.5 -.866025)
\move(2.5 -7.794225)
\rlvec(.5 -.866025) 
\move(1.5 -6.062175)
\rlvec(1 0)
\rmove(.5 -.866025)
\rlvec(1 0) \rmove(.5 -.866025)
\rlvec(-.5 -.866025) 
\htext(-5.5 -11.4){\small a. The perfect matching corresponding}
\htext(-5.5 -12){\small \phantom{a. } to the rhombus tiling of
Figure~\ref{fliegefi}. }
\move(12 0)
\bsegment
\drawdim cm
\setunitscale.5773503
\linewd.04
\hex \hex \hex
\move(-1.5 -.866025)
\hex \hex \hex
\move(-3 -1.73205)
\hex \hex 
\move(-3 -3.4641)
\rlvec(-.5 -.866025)

\move(-3 -5.19615)
\hex \hex \hex \hex
\move(-3 -6.9282)
\hex \hex \hex
\move(1.5 -6.062168)
\hex \hex
\move(4 -3.4641)
\rlvec(.5 -.866025)
\rmove(-.5 -.866025)
\rlvec(.5 -.866025)
\rlvec(-.5 -.866025)
\htext(-2 -11.4){\small b. $G^+$ and $G^-$.}
\esegment
}
\caption{\small \label{hexfi} }
\end{figure}

The number of lattice paths with positive horizontal and
negative vertical steps from $(a,b)$ to $(c,d)$ equals 
$\binom{c-a+b-d}{b-d}$.
Therefore, the number of families of nonintersecting 
lattice paths (equivalently, the number of rhombus tilings of $R^+$) 
is equal to the following determinant:

$$
\det_{1\le i,j\le n}(\mathcal{P}(P_i\rightarrow Q_j))=\det_{1\le i,j\le n}
\(\binom{m+j-1}{m-j+i}\).
$$
This proves Lemma~\ref{plusdet}.
\end{proof}
\begin{figure}
\centertexdraw{
\drawdim truecm \linewd.02
\rhombus \rhombus \rhombus \ldreieck
\move(-.866025 -.5) 
\rhombus \rhombus \rhombus \rhombus \ldreieck
\move(-1.73205 -1) 
\rhombus \rhombus \rhombus \rhombus 
\move(-1.73205 -1)
\rdreieck \rhombus \rhombus
\move(-1.73205 -2)
\rdreieck 
\htext(-2 -4){\small a. $R^+$, the upper half of the hexagon.}
\move(6 0)
\bsegment
\linewd.02
\rhombus \rhombus \rhombus \rhombus \rhombus
\move(0 0)
\rdreieck \rhombus \rhombus \rhombus \rhombus
\move(0 -1)
\rdreieck \rhombus \rhombus \rhombus
\move(2.598075 -.5)
\rhombus \rhombus \ldreieck
\move(1.73205 0) \ldreieck
\move(4.3301 0.5) \rdreieck \ldreieck
\esegment
\htext(6 -4){\small b. $R^-$, the lower half of the hexagon.}
}
\caption{\label{uplowfi}}
\end{figure}
\begin{figure}
\vbox{
\centertexdraw{
\drawdim truecm 
\linewd.05
\move(0 0)
\RhombusA \RhombusB \RhombusA \RhombusB
\move(-.866025 -.5)
\RhombusA \RhombusB \RhombusB
\move(-1.73205 -1)
\RhombusB \RhombusB
\move(1.73205 0) \RhombusC \RhombusC
\move(2.598075 -1.5) \RhombusC
\move(-.866025 -.5) \RhombusC
\move(-.866025 -1.5) \RhombusC
\move(.866025 -1.5) \RhombusC
\ringerl(.433012 .25) 
\hdSchritt \vdSchritt \hdSchritt \vdSchritt 
\ringerl(-.433012 -.25)
\hdSchritt \vdSchritt \vdSchritt 
\ringerl(-1.299037 -.75)
\vdSchritt \vdSchritt 
\ringerl(2.165062 -2.75) \ringerl(.433012 -2.75) \ringerl(-1.299038
-2.75)
\htext(-2 -4){\small a. A tiling of the upper half 
of the hexagon}
\htext(-2 -4.5){\small and the corresponding lattice path family.}
\move(9 0)
\bsegment
\ringerl(.433012 .25) 
\hdSchritt \vdSchritt \hdSchritt \vdSchritt
\ringerl(-.433012 -.25)
\hdSchritt \vdSchritt \vdSchritt 
\ringerl(-1.299037 -.75)
\vdSchritt \vdSchritt
\ringerl(2.165062 -2.75) \ringerl(.433012 -2.75) \ringerl(-1.299038
-2.75)
\esegment
\htext(7.5 -4){\small b. The paths isolated.}
\htext(3 -13){\small c. The corresponding lattice path family.}}
\vskip-8cm
$$
\Einheit=1cm
\Gitter(5,5)(0,0)
\Koordinatenachsen(5,5)(0,0)
\Pfad(0,0),22\endPfad
\Pfad(2,1),22\endPfad
\Pfad(1,3),1\endPfad
\Pfad(2,4),1\endPfad
\Pfad(3,3),1\endPfad
\Pfad(3,3),2\endPfad
\Pfad(4,2),2\endPfad
\Kreis(0,0.02) \Kreis(0,2.02) \Kreis(2,1.02)
\Kreis(1,3.02) \Kreis(4,2.02)
\Kreis(2,4.02)
\Label\lo{P_1}(0,2)
\Label\lo{P_2}(1,3)
\Label\lo{P_3}(2,4)
\Label\ru{Q_1}(0,0)
\Label\ru{Q_2}(2,1)
\Label\ru{Q_3}(4,2)
\hskip2cm
$$
\vskip2cm}
\caption{\label{upfi}}
\end{figure}

This determinant can be evaluated with the help of
the following determinant lemma (\cite{kratt}, Lemma~2.2).
\begin{lemma} \label{det}
\begin {multline*} 
\det_{1 \le i,j \le n} \big( (x_j+a_n)(x_j+a_{n-1}) \cdots
(x_j+a_{i+1})(x_j+b_i)(x_j+b_{i-1}) \dots (x_j+b_2) \big) \\
=\prod _{1 \le i < j \le n} ^{} {(x_i-x_j)} \prod _{2 \le i \le j
\le n} ^{}{(b_i-a_j)}.
\end{multline*}
\end{lemma}
The determinant in Lemma~\ref{plusdet} factors as follows.
\begin{lemma} \label{g+}
\begin{equation*}M(G^+)= M(R^+)=\frac{\h(n)\prod _{2\le i \le j \le n} ^{}{(2m+2j-i)}}
{\prod _{j=1} ^{n}{(2j-2)!}}.
\end{equation*}
\end{lemma}
Here, in abuse of notation, $M(R)$ denotes the number of rhombus tilings 
of $R$, if $R$ is a region of triangles.
\begin{proof}
Before we can apply Lemma~\ref{det}, we have to transform 
the expression of Lemma~\ref{plusdet} in the following way.
\begin{multline*} 
\begin{split}
M(G^+)&=\det_{1\le i,j\le n}\(\binom{m+j-1}{m-j+i}\)\\
&=\prod _{i=1} ^{n}{(-2)^{i-1}}
\prod _{j=1} ^{n}{\frac {(m+j-1)!} {(n+m-j)!(2j-2)!}}\end{split}
\\
\times\det_{1\le i,j\le n} 
\Big((m+i+1-j)(m+i+2-j)\dots(m+n-j) \\
 \left.\cdot\(-j+\frac {i}{2}\)\(-j+\frac {i}{2}-\frac
{1}{2}\)
\dots\(-j+1\)\right).
\end{multline*} 
Now we apply Lemma~\ref{det} with $x_k=-k,\,a_k=m+k,\,b_k=\frac {k} {2}$
and simplify to get the claimed result.
\end{proof}
\end{section}
\begin{section}{The matchings count for the lower half}\label{lowsec}
\begin{figure}
\vbox{
\centertexdraw{
\drawdim truecm 
\linewd.05
\move(0 0)
\RhombusA \RhombusB \RhombusB \RhombusA \RhombusA 
\move(1.72305 0)
\RhombusB \RhombusA \RhombusB \RhombusA 
\move(2.598075 -.5)
\RhombusA \RhombusA \RhombusB 
\move(0 0)
\RhombusC
\move(0 -1) \RhombusC
\move(4.33012 .5)
\RhombusC
\ringerl(.433012 .25) 
\hdSchritt \vdSchritt \vdSchritt \hdSchritt \hdSchritt \fcir f:0 r:.07
\ringerl(2.165063 .25)
\vdSchritt \hdSchritt \vdSchritt \hdSchritt \fcir f:0 r:.07
\ringerl(3.032088 -.25)
\hdSchritt \hdSchritt \vdSchritt \fcir f:0 r:.07
\htext(-.5 -4){\small a. A tiling of the lower half 
of the hexagon}
\htext(-.5 -4.5){\small and the corresponding lattice path family.}
\move(9 1)
\bsegment
\ringerl(-.433012 -.5) 
\hdSchritt\vdSchritt\vdSchritt\hdSchritt\hdSchritt \fcir f:0 r:.07
\ringerl(1.298038 -.5)
\vdSchritt\hdSchritt\vdSchritt\hdSchritt\fcir f:0 r:.07
\ringerl(2.165063 -1)
\hdSchritt\hdSchritt\vdSchritt\fcir f:0 r:.07
\esegment
\htext(9.5 -4){\small b. The paths isolated.}
\htext(4.5 -12){\small c. The corresponding lattice path family.}}
\vskip-6.5cm
$$
\Einheit=1cm
\Gitter(6,4)(0,0)
\Koordinatenachsen(6,4)(0,0)
\Pfad(0,2),1\endPfad
\Pfad(1,0),22\endPfad
\Pfad(1,0),11\endPfad
\Pfad(2,2),2\endPfad
\Pfad(2,2),1\endPfad
\Pfad(3,1),2\endPfad
\Pfad(3,3),11\endPfad
\Pfad(5,2),2\endPfad
\Pfad(3,1),1\endPfad
\Kreis(0,2.02) \Kreis(3,0.02) \Kreis(4,1.02)
\Kreis(2,3.02) \Kreis(3,3.02)
\Kreis(5,2.02) 
\Label\lo{R_1}(0,2)
\Label\lo{R_2}(2,3)
\Label\lo{R_3}(3,3)
\Label\ru{S_1}(3,0)
\Label\ru{S_2}(4,1)
\Label\ru{S_3}(5,2)
\Label\o{\frac {1} {2}}(.5,2.1)
\hskip2cm
$$
\vskip2cm}
\caption{\label{lowfi}}
\end{figure}
In this section we reduce the computation of $M(G^-)$ to the
evaluation of a determinant, which is stated in the following lemma.
\begin{lemma}\label{minusdet}
\begin{align*} 
M(G^-)=M(R^-)
=\frac{(n+m-s)(s+m)} {(2n-2s)\prod _{i=1} ^{n}{(2n+1-2i)!}} 
\times\det (B_{ij}),\quad s=0,\dots,n-1
\end{align*}
where  $R^-$ is the region of triangles which is described in the
proof of the lemma, and where
\begin{equation} \label{b1}
B_{ij}=\begin{cases}
(n+2+j-2i)_{n-j}(i+m+1-j)_{j-1}\left(m+\frac {n}{2}+\frac {1}{2}-\frac
{j}{2}\right)
\quad & i\not=s+1 \\
(n+1+j-2s)_{n-j}(s+m+1-j)_{j-1} &i=s+1.
\end{cases}
\end{equation} 
\end{lemma}

Here $(a)_n:=a(a+1)(a+2)\dots (a+n-1)$ is the usual shifted
factorial. The entry for $i\not=s+1$ can also be written as 
\begin{equation}\label{balt}
B_{ij}=\frac {1} {2}(n+1+j-2i)_{n-j+1}(i+m+1-j)_{j-1}
+(n+2+j-2i)_{n-j}(i+m-j)_{j}.
\end{equation}

\begin{proof}
We start analogously to Section~\ref{upsec} and convert $G^-$
(exemplified in Figure~\ref{hexfi}b) back
to a region $R^-$ of triangles (see Figure~\ref{uplowfi}b), so that
the perfect matchings of $G^-$ correspond bijectively to the rhombus
tilings of $R^-$.
However, since $G^-$ contains edges on the symmetry axis of $G$,
which, by Lemma~\ref{fact}, has the consequence that they count with
weight $\frac {1} {2}$ in $G^-$, we are dealing with a weighted count
of the rhombus tilings of $R^-$, where rhombi such as the top-left
rhombus in Figure~\ref{lowfi}a count with weight $\frac {1} {2}$.
Again, we count the rhombus tilings of $R^-$ by counting
the number of nonintersecting lattice path families leading from
upper left to lower right edges. The starting and
end points can be easily read off Figure~\ref{lowfi}.

The starting points are (see Figure~\ref{lowfi} for the meaning of $s$, 
$1\le s \le n-1$):
$$R_i=
\begin{cases}
(2i-2,m+i-1)\quad &\text {for } i\not=s+1\\  
(2s-1,m+s-1)&\text {for }i=s+1.
\end{cases}
$$

The end points are:
$$S_j=(n+j-1,j-1)\quad \text {for }j=1,\dots,n.$$

Now we apply again the main result for nonintersecting lattice paths.
The matrix entries are $\mathcal{P}(R_i\rightarrow S_j)$.
We note that a positive horizontal step starting at $R_i$, 
$i\not=s+1$ corresponds to a rhombus of weight $\frac {1} {2}$ in
$R^-$, so paths starting with a horizontal step at $R_i$, $i\not =
s+1$ are counted with weight $\frac {1} {2}$.
Therefore, we count paths starting with a horizontal step and 
paths starting with a vertical step separately. 

We get $M(G^-)=\det_{1\le i,j \le n} \( A_{ij} \)$, where
\begin{align} \label{aeq} A_{ij}&=
\mathcal{P}(R_i\rightarrow S_j)=
\begin{cases}
\dfrac {1} {2}\dbinom{n+m-i}{m+i-j}+\dbinom{n+m-i}{m+i-1-j}\quad &\text
{for }i\not=s+1\\
\dbinom{n+m-s}{m+s-j}&\text {for }i=s+1.
\end{cases}
\end{align}
Since this expression also makes sense for $s=0$, we can include 
this case in the following calculations. 
We pull factors out of the rows of $\det \(A_{ij}\)$, 
so that the remaining entries are polynomials in $m$ and get the
expression in Lemma~\ref{minusdet}.
\end{proof}

The next lemma states an evaluation of the determinant in
Lemma~\ref{minusdet}. An outline of the proof is provided immediately
below. The details are given in Sections~\ref{halbsec} -- \ref{degsec}.
\begin{lemma} \label{g-}
\begin{multline*} 
\det\(B_{ij}\)=\frac {2^{\binom{n-1}{2}}\h(n)(2n-2s-1)!!
(2s-1)!!}{(n-s-1)! s! (m+s)(m+n-s)}\\
\hskip1.5cm\times
\prod _{k=1} ^{n-2}{\left( m+k+\frac{1}{2}\right)^{\min(k,n-1-k)} 
\prod_{k=0} ^{n}{(m+k)^{\min(k+1,n-k+1)}}}.
\end{multline*}
\end{lemma}

{\em Outline of the proof of Lemma~\ref{g-}:}
In Section~\ref{halbsec} we prove that $\prod _{k=1} ^{n-2}
{\left( m+k+\frac{1}{2}\right)^{\min(k,n-1-k)}}$ divides $\det
\(B_{ij}\)$ as a polynomial in $m$.
In Section~\ref{ganzsec} we prove that $\prod_{k=0} ^{n}
{(m+k)^{\min(k+1,n-k+1)}}$ divides $\det \(B_{ij}\)$ as a 
polynomial in $m$.

The results are stated in Lemmas~\ref{halb} and
\ref{ganz}. We show this with the help of linear combinations of
rows and columns, which vanish, if one of the linear factors is set
equal to zero. 

In Section~\ref{degsec} we compute the degree of the determinant as a
polynomial in $m$. It is exactly equal to the number of linear
factors we have already found to divide the determinant. So we know
the determinant up to a constant factor. 
We compute this constant in Lemma~\ref{lead} 
by replacing each entry by its leading coefficient and using
Vandermonde's determinant formula.

Lemmas~\ref{halb}, \ref{ganz} and \ref{lead} immediately give
Lemma~\ref{g-}. \qed

\end{section}

\begin{section}{The ``half-integral'' factors of $\det\(B_{ij}\)$} 
\label{halbsec}
In this section we prove the following lemma (see Lemma~\ref{minusdet}
for the definition of $B_{ij}$).
\begin{lemma} \label{halb}
$\prod _{k=1} ^{n-2}{\left( m+k+\frac{1}{2}\right)^{\min(k,n-1-k)}}$
divides $\det \( B_{ij}\)$ as a polynomial in $m$. 
\end{lemma}
\begin{proof}
We find linear combinations of columns which give zero for
$m=-\left(k+\frac{1}{2}\right)$. First, we show that the
following linear combination of columns 
equals zero for $i\not=s+1$, $0\le l\le k$, $l\ge 2k-n+1$:
\begin{equation} \label{halbsum}
\sum _{j=0} ^{l}{\binom{l}{j}
B_{i,n+2l-2k-j}\Big|_{m=-k-\frac {1} {2}} }=0.
\end{equation}

In order to establish this, we break the sum in two parts according 
to the two summands of $B_{ij}$ in equation~\eqref{balt}
and convert them to hypergeometric form.
The left-hand side of \eqref{halbsum} becomes
\begin{multline*} 
({k-2l+i-1-\frac {1} {2}-n})_{n-2k+2l}
(2-2k+2l-2i+2n)_{2k-2l}\,
{}_{2}F_{1}\!\left[\begin{matrix}{-l,-1+2k-2l+2i-2n}\\{
2k-2l+i-k-\frac {1} {2}-n}\end{matrix};{\displaystyle1}\right]\\
+
({1+2k-2l+i-k-\frac {1} {2}-n})_{-1-2k+2l+n}
({2-2k+2l-2i+2n})_{2k-2l}\\
\times
{}_{2}F_{1}\!\left[\begin{matrix}{-l,-1+2k-2l+2i-2n}\\{1+
2k-2l+i-k-\frac {1} {2}-n}\end{matrix};{\displaystyle1}\right]. 
\end{multline*}

Now, Vandermonde's summation formula 
$${} _{2} F _{1} \!\left [ \begin{matrix} { a, -n}\\ { c}\end{matrix} ;
  {\displaystyle 1}\right ] = 
 {\frac{({ \textstyle -a + c}) _{n} }{({ \textstyle c}) _{n} }}$$
is applicable to both ${} _{2} F _{1}$--series, since $l \ge 0$. 
It is directly verifiable that the two resulting expressions sum to zero.

It is easily seen that the conditions on $l$ in \eqref{halbsum} allow
$\min(k+1,n-k)$ possible values for $l$.
Thus, we have $\min(k+1,n-k)$ independent 
linear combinations of columns which vanish in all  coordinates 
except possibly in the $(s+1)$-th coordinate. 
(Recall that \eqref{halbsum} is valid only for
$i\not= s+1$). It is clear that an appropriate combination of two of 
these linear combinations vanishes in every coordinate. 
So we have $\min(k,n-k-1)$ independent linear combinations vanishing
at $m=-k-\frac {1} {2}$ for $k=1,\dots,n-2$, which proves
Lemma~\ref{halb}.  
\end{proof}
\end{section}
\begin{section}{The ``integral'' factors of $\det\(B_{ij}\)$}
\label{ganzsec}
In this section we prove the following lemma (see Lemma~\ref{minusdet}
for the definition of $B_{ij}$).
\begin{lemma}\label{ganz}
$\frac{1}
{(n+m-s)(s+m)} \prod_{k=0} ^{n}{(m+k)^{\min(k+1,n-k+1)}}$ divides 
$\det \(B_{ij} \)$ as a polynomial in $m$.
\end{lemma}
\begin{proof}
We use linear combinations of the rows of $B_{ij}$ that vanish for
$m=-k$. 
Without loss of generality, we can assume that $s\le \frac {n} {2}$ 
because both the final result (see Theorem~\ref{th1}) 
and the number of rhombus tilings of the original
graph are invariant under the transformation $s\rightarrow n-s$.

Most of the factors $(m+k)$ can be pulled out directly from the rows of 
$B_{ij}$. 
In fact, it is easily seen that row $i$ is divisible by 
$(m+1-i+n)_{2i-n-1}$ for $2i\ge n+2$. 
The product of these terms equals
$\prod _{k=1} ^{n-1}{(m+k)^{\min(k,n-k)}}$.
The matrix $\(C_{ij}\)$ which remains after pulling out these factors
from $B_{ij}$ looks as follows:
$$
C_{ij}=\begin{cases}
(n+j+1-2s)_{n-j}(s+m+1-j)_{j-1} &i=s+1,\\
(n+2+j-2i)_{n-j}(i+m+1-j)_{j-2i+n}(2m+n+1-j) &i\not=s+1,2i\ge n+2,\\
(n+2+j-2i)_{n-j}(i+m+1-j)_{j-1}(2m+n+1-j) &i\not=s+1,2i< n+2.
\end{cases}
$$
To finish the proof of Lemma~\ref{ganz} we have to find for each
value $k=0,\dots,n$, $k\not=s$ and $k\not=(n-s)$ one vanishing linear
combination of the rows of $C_{ij}$.

We start with the case $k< s$. 
We claim that 
\begin{multline} \label{linganz1}
\sum_{i=k+1}^{s}
{\left((-1)^{i-k+1}\binom{s-k-1}{i-k-1}\frac {(n+\frac {3}{2}-i)_{i-k-1}
(n-i+1)_{i-k-1} }{(s+\frac {1}{2}-i)_{i-k-1}
(n-k-i+1)_{i-k-1}}  C_{ij} \Big|_{m=-k} \right)} \\
\hskip-2cm+
(-1)^{s-k+2}
\frac {2(n+\frac {3}{2}-s-1)_{s-k}(n-s+1)_{s-k-1}}
{(\frac {1}{2})_{s-k-1} (n-k-s+1)_{s-k-1}} C_{s+1,j}\Big|_{m=-k}=0.
\end{multline}
If $-k+s-j<0$ then the terms $(i-k+1-j)_{j-1}$ and $(s-k+1-j)_{j-1}$
(which are factors of $C_{ij}\Big|_{m=-k}$ and $
C_{s+1,j}\Big|_{m=-k}$, respectively) are zero for all 
occurring indices.
If $-k+s-j\ge 0$, we reverse the order of summation in the sum and
write \eqref{linganz1} in hypergeometric form,
\begin{multline*}
\left( -1 \right)^{-k + s} \frac { 
   (n + \frac{3}{2} - s) _{s-k -1} \,
   (1-j-k+s) _{j-1} } 
  {(\frac{1}{2}) _{s-k-1} \,
(n-k-s+1)_{s-k-1} ( 2 + j + n - 2\,s) _{n-j} \,
   (n-s+1) _{s-k-1} }\\
\times (-1+j+2k-n) \,
   {} _{3} F _{2}\!\left[\begin{matrix}{\frac{1}{2},
n-k+1-s,j+k-s}\\{1+\frac{j}{2}+\frac{n}{2}-s,\frac{3}{2}+
\frac{j}{2}+
\frac{n}{2}-s}\end{matrix};{\displaystyle1}\right]
\\\hskip-2cm+(-1)^{s-k+2}
\frac{2(n+\frac{3}{2}-s-1)_{s-k}(n-s+1)_{s-k-1}}
{(\frac{1}{2})_{s-k-1}(n-k-s+1)_{s-k-1}} C_{s+1,j}\Big|_{m=-k}=0.
\end{multline*}

Now we can apply the Pfaff--Saalsch\"utz summation formula 
(\cite{SlatAC}, (2.3.1.3); Appendix (III.2)),
\begin{equation}\label{pfaff}
{} _{3} F _{2} \!\left [ \begin{matrix} { a, b, -n}\\ { c, 1 + a + b - c -
  n}\end{matrix} ; {\displaystyle 1}\right ] = 
 {\frac{({ \textstyle -a + c}) _{n} \,({ \textstyle -b + c}) _{n} } 
  {({ \textstyle c}) _{n} \,({ \textstyle -a - b + c}) _{n} }}.
\end{equation}
It is easily verified that the 
resulting sum of two terms equals zero.

The case $k> n-s$ is quite similar. 
We claim that
\begin{multline} \label{linganz2} 
\sum_{i=n-k+1}^{s}
{\left((-1)^{i-n+k+1}\binom{s-n+k-1}{i-n+k-1}\frac {(n+\frac
{3}{2}-i)_{i-n+k-1}
(n-i+1)_{i-n+k-1} }{(s+\frac {1}{2}-i)_{i-n+k-1}
(k-i+1)_{i-n+k-1}} C_{ij}\Big|_{m=-k}\right)} \\
\\\hskip-2cm-
(-1)^{s-n+k+2}
 \frac {2(n+\frac {3}{2}-s-1)_{s-n+k}(n-s+1)_{s-n+k-1}}
{(\frac {1}{2})_{s-n+k-1} (k-s+1)_{s-n+k-1}} 
C_{s+1,j}\Big|_{m=-k}=0.
\end{multline}
Converting the reversed sum to hypergeometric form gives
\begin{multline*}
(-1)^{s+k-n}(-1+j+2k-n)(2+j+n-2\,s)_{n-j}\,
(1+n-s)_{-1+k-n+s}\\
\times\frac{
(\frac{3}{2}+n-s)_{-1+k-n+s}
(1-j-k+s)_{-1+j}
}
{
(\frac{1}{2})_{-1+k-n+s}(1+k-s)_{-1+k-n+s}
}\\
\times{}_{3}F_{2}\!
\left[
\begin{matrix}{1-k+n-s,\frac{1}{2},j+k-s}\\
{1+\frac{j}{2}+\frac{n}{2}-s,\frac{3}{2}+\frac{j}{2}+\frac{n}{2}-s}
\end{matrix};
{\displaystyle1}\right].
\end{multline*}
Again, the Pfaff--Saalsch\"utz summation formula \eqref{pfaff} 
is applicable because $-(1-k+n-s)$ is a nonnegative integer. It 
is easily checked that the resulting terms sum to 0. 

So our remaining task is the case $s<k<n-s$. 
For $s<k\le \frac {n} {2}$ we consider the following linear
combination,

\begin{multline} \label{linganz3}
\sum_{i=k+1}^{\left\lfloor\frac {n+1}{2}\right\rfloor} 
{(-4)^{n-i}
\frac {(s-i+1)_{i-k-1}(s-n+\frac {1} {2})_{n-k-1}}
{(2n-2i+1)!(s+\frac {1}{2}-i)_{i-k-1}
(n+1-s)_{-k}
}
(i-k)_{n+1-2i} C_{ij}\Big|_{m=-k} 
}\\
+
\sum_{i=\left\lfloor\frac{n+3}{2}\right\rfloor}
^{\left\lfloor\frac {n+1+j}{2}\right\rfloor} 
{(-4)^{n-i}
\frac {(s-i+1)_{i-k-1}(s-n+\frac {1} {2})_{n-k-1}}
{(2n-2i+1)!(s+\frac {1}{2}-i)_{i-k-1}
(n+1-s)_{-k}
}
C_{ij} \Big|_{m=-k}
}\\
-C_{s+1,j}\Big|_{m=-k} =0.
\end{multline} 

Now, both the term $(i-k)_{n+1-2i} C_{ij}\Big|_{m=-k}$ which is
part of the first sum in \eqref{linganz3} and the term
$C_{ij} \Big|_{m=-k}$ which is part of the second sum in
\eqref{linganz3} are equal to
$(n+2+j-2i)_{n-j}(i-k+1-j)_{j-2i+n}(-2k+n+1-j)$, so we can combine
the two sums into one sum of a hypergeometric term.

We distinguish two cases according to the
parity of $n-j$. In both cases we reverse the sum and convert it to
hypergeometric form.
The resulting two hypergeometric series are
\begin{align*}
{_3}F_2
&\left[
\begin{matrix}1+l+m,-l-m,\frac{1}{2}+l-n+s\\
\frac{3}{2},1+l-n+s
\end{matrix}
;1
\right]&\text{for }n-j&=2l,\\
\ _{3} ^{}
F_2&\left[\begin{matrix}1+l+m,-1-l-m,
\frac{1}{2}+l-n+s\\
\frac{1}{2},1+l-n+s\end{matrix};1\right]\quad
&\text{for }n-j&=2l+1.
\end{align*}
So we can use the Pfaff--Saalsch\"utz summation formula \eqref{pfaff}
again.
It is easily verified that in both cases the resulting terms add to
zero.

The case $\frac {n} {2}<k<n-s$ is handled similarly. We claim that
the following sum equals zero,

\begin{multline}\label{linganz4}
\sum _{i=n-k+1} ^{\left\lfloor\frac {n+1} {2}\right\rfloor}
{(-4)^{n-i}\frac
{(s-i+1)_{i-n+k-1}(s-n+\frac {1} {2})_{k-1}}
{(2n-2i+1)!(s+\frac{1}{2}-i)_{i-n+k-1}
(n+1-s)_{-n+k}
}
 (i-k)_{n+1-2i}  C_{ij}\Big|_{m=-k}}\\
+
\sum _{i=\left\lfloor\frac {n+3} {2}\right\rfloor} 
^{\left\lfloor\frac {n+j+1} {2}\right\rfloor}
{(-4)^{n-i}\frac
{(s-i+1)_{i-n+k-1}(s-n+\frac {1} {2})_{k-1}}
{(2n-2i+1)!(s+\frac{1}{2}-i)_{i-n+k-1}
(n+1-s)_{-n+k}
}C_{ij}\Big|_{m=-k}}\\
- (-1)^n C_{s+1,j}\Big|_{m=-k} =0.
\end{multline}
Again, we write the two sums as one single sum, 
distinguish two cases according to the parity of $n-j$, and reverse 
the order of summation. Conversion to hypergeometric form of the
resulting sums gives
\begin{align*}
 _{3} ^{}F_2&\left[\begin{matrix}1+l-k,-l+k,
\frac{1}{2}+l-n+s\\
\frac{3}{2},1+l-n+s\end{matrix};1\right]
&\text{for }n-j&=2l,\\
 _{3} ^{}F_2&\left[\begin{matrix}1+l-k,-1-l+k,
\frac{1}{2}+l-n+s\\
\frac{1}{2},1+l-n+s\end{matrix};1\right]
&\text{for }n-j&=2l+1.
\end{align*}
The Pfaff--Saalsch\"utz summation formula \eqref{pfaff}
can be applied in both cases. It is easily seen that the results
vanish after subtraction of $(-1)^n C_{s+1,j}\Big|_{m=-k} $.
Thus Lemma~\ref{ganz} is proved.
\end{proof}
\end{section}

\begin{section}{The degree and the leading coefficient}\label{degsec}
We have to find the degree and the leading coefficient of the determinant
$\det \( B_{ij}\)$ as a polynomial in $m$. The degree of $B_{ij}$ is 
$j-1$ for $i=s+1$ and $j$ else. The degree of $\det\(B_{ij}\)$ is clearly
$\binom{n+1}{2}-1,$ which is easily seen to be the number of linear
factors we have found to divide $\det\(A_{ij}\)$.
Therefore $\det\(B_{ij}\)$ is the product of the linear factors and
the leading coefficient.

To compute the leading coefficient we look at the leading coefficient of 
each entry. The leading coefficients of each entry give the matrix $D$, 
with
$$
D_{ij}=(x_i+n+j)_{n-j} \text {, where } 
x_i=\left\{ \begin{matrix} 2-2i&i\not=s+1\\1-2s&i=s+1.\end{matrix} \right.
$$
This matrix can be transformed by column reduction to $(x_i^{n-j})$,
but this is just the Vandermonde determinant $\prod
_{1\le i<j\le n} ^{}{(x_i-x_j)}$. Plugging in the values of the $x_i$
gives the following lemma.
\begin{lemma} \label{lead}
The leading coefficient of $\det \(B_{ij}\)$ is
$$
\frac {2^{\binom{n-1}{2}}\h(n)(2n-2s-1)!!
(2s-1)!!}{(n-s-1)! s!}.
$$
\end{lemma}

Lemmas~\ref{halb}, \ref{ganz} and \ref{lead} immediately give
Lemma~\ref{g-}. Thus the proof of Lemma~\ref{g-} is complete. \qed

\end{section}
\begin{section}{Proof of Theorem~\ref{th2}}\label{oddsec}
\begin{figure}
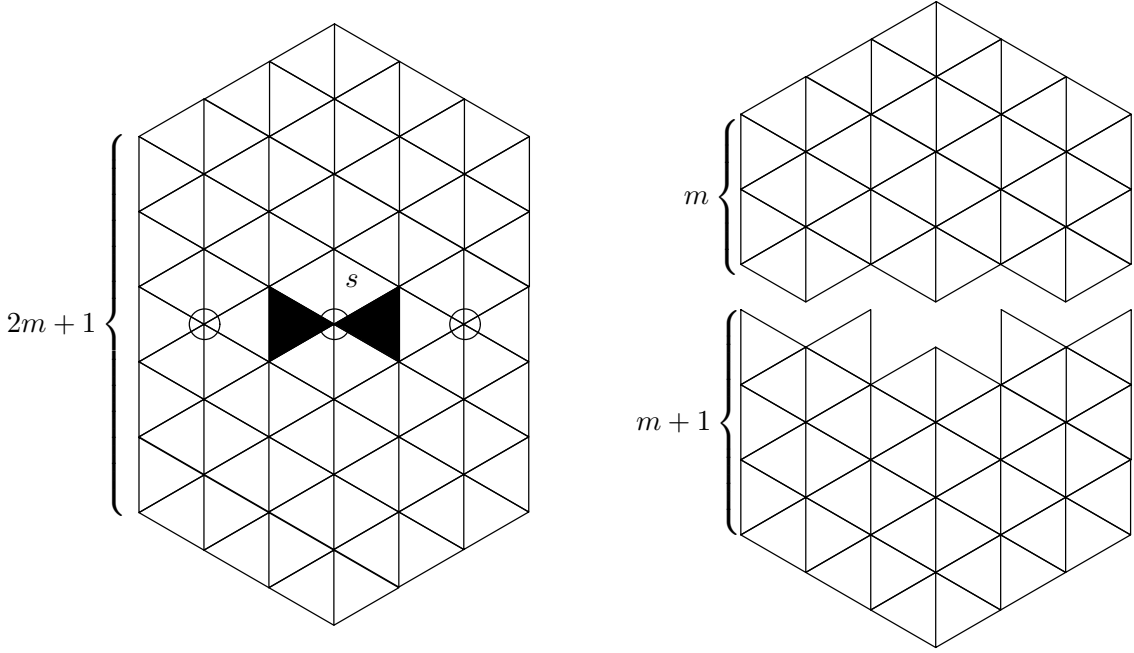

\centertexdraw{
\drawdim truecm \linewd.02
\rhombus \rhombus \rhombus \ldreieck 
\move(-.866025 -.5)
\rhombus \rhombus \rhombus \rhombus \ldreieck
\move(-1.732 -1)
\rhombus \rhombus \rhombus \rhombus \rhombus \ldreieck
\move(-1.732 -1) \rdreieck \rhombus \rhombus \rhombus \rhombus \rhombus
\ldreieck
\move(-1.732 -2) \rdreieck \rhombus \rhombus \rhombus \rhombus
\rhombus \ldreieck 
\move(-1.732 -3)	
\rdreieck \rhombus \rhombus \rhombus \rhombus \rhombus 
\move(-1.732 -4)
\rdreieck \rhombus \rhombus \rhombus \rhombus
\move(-1.732 -5)
\rdreieck \rhombus \rhombus \rhombus
\move(.866025 -3.5)
\bsegment
\rlvec(.866 .5) \rlvec(0 -1) \lfill f:0
\move(0 0)
\rlvec(-.866 .5) \rlvec(0 -1) \lfill f:0
\esegment
\move(-.866025 -3.5)
\lcir r:.2
\move(.866025 -3.5)
\lcir r:.2
\move(2.598075 -3.5)
\lcir r:.2

\rtext td:0 (1 -3){$s$}
\rtext td:0 (-3.5 -3.6) {$2m+1 \left\{ \vbox{\vskip2.7cm}\right.$}

\move(8 0.3)
\bsegment
\rhombus \rhombus \rhombus \ldreieck
\move(-.866025 -.5) 
\rhombus \rhombus \rhombus \rhombus \ldreieck
\move(-1.73205 -1) 
\rhombus \rhombus \rhombus \rhombus \rhombus
\move(-1.73205 -1)
\rdreieck \rhombus \rhombus \rhombus
\move(-1.73205 -2)
\rdreieck \rhombus
\rtext td:0 (-2.5 -2.2) {$m \left\{ \vbox{\vskip1.2cm}\right.$}
\move(-1.73205 -3.6)
\bsegment
\rdreieck \rhombus \rhombus \rhombus \rhombus \rhombus 
\move(0 -1)
\rdreieck \rhombus \rhombus \rhombus \rhombus
\move(0 -2)
\rdreieck \rhombus \rhombus \rhombus
\move(.866025 -.5) \ldreieck \rhombus \rhombus \rhombus \ldreieck
\move(3.463 0) \rdreieck \rhombus \ldreieck
\rmove(-.866025 1.5) \ldreieck
\rtext td:0 (-1.4 -1.6) {$m+1 \left\{ \vbox{\vskip1.6cm}\right.$}
\esegment
\esegment
} 
\caption{\label{oddfi} The hexagon in the case of odd sidelength
and the two halves $R^+$ and $R^-$. $m=2$, $n=3$, $s=2$.}
\end{figure}
If the side divided by the symmetry axis has odd length, the position
$s$ of the missing triangles ranges from $1$ to $n$ (see
Figure~\ref{oddfi}). 
We can form the inner dual graph and denote it by $\odd{G}$. Now we
can proceed analogously to Section~\ref{cutsec}, 
break the graph in two parts $\odd{G^+}$ and
$\odd{G^-}$ with the help of the matchings factorization theorem (see
Lemma~\ref{fact}). We convert $\odd{G^+}$ and
$\odd{G^-}$ back to regions $\odd{R^+}$ and $\odd{R^-}$ of triangles
and have to count rhombus tilings again. $\odd{R^+}$ and $\odd{R^-}$
are shown in an example in Figure~\ref{oddfi}. 

Thus, we have 
\begin{equation} \label{facto}
\odd{G}=2^{n-1}M(\odd{R^+})M(\odd{R^-}).
\end{equation} 

Now we reduce the evaluation of $M(\odd{R^+})$ to the evaluation
of $M(R^+)$, which we have already done in Lemma~\ref{g+}.
$M(R^+)$ and $M(\odd{R^+})$ are related in the following way.
The tiles of the upper half $R^+$ of the hexagon with sides $n,n,2m,n,n,2m$
(as exemplified in Figure~\ref{fliegefi}) 
sharing an edge with the border of length $m$ are 
enforced as shown in Figure~\ref{relupfi} (the forced tiles are shaded).
After removal of these tiles we are left with the upper half
$\odd{R^+}$ of the hexagon with sides $n-1,n-1,2m+1,n-1,n-1,2m+1$.
\begin{figure}
\centertexdraw{
\drawdim truecm \linewd.02
\rhombus \rhombus \rhombus \ldreieck
\move(-.866025 -.5) 
\rhombus \rhombus \rhombus \rhombus \ldreieck
\move(-1.73205 -1) 
\rhombus \rhombus \rhombus \rhombus \rhombus 
\move(-1.73205 -1)
\rdreieck \rhombus \rhombus \rhombus
\move(-1.73205 -2)
\rdreieck \rhombus
\move(-1.73205 -1)
\rlvec(0 -1) \rlvec(-.866025 -.5) \lfill f:.5
\rlvec(.866025 .5) \rlvec(0 -1) \lfill f:.5
\rlvec(-.866025 -.5) \rlvec(0 1) \lfill f:.5
\rlvec(0 1) \rlvec(.866025 .5)  \lfill f:.5
\move(3.465 -1)
\rlvec(0 -1) \rlvec(.866025 -.5) \lfill f:.5
\rlvec(-.866025 .5) \rlvec(0 -1) \lfill f:.5
\rlvec(.866025 -.5) \rlvec(0 1) \lfill f:.5
\rlvec(0 1) \rlvec(-.866025 .5)  \lfill f:.5 }
\caption{\label{relupfi} $M(\odd{R^+}(3,2))$ equals $M(R^+(2,2))$.}
\end{figure}
Thus, we have $M(R^+(n,m))=M(\odd{R^+}(n-1,m))$ and Lemma~\ref{g+} 
implies directly the following result:
\begin{equation}\label{odd+}
M(\odd{R^+})=
\frac{\h(n+1)\prod _{2\le i \le j \le n+1} ^{}{(2m+2j-i)}}
{\prod _{j=1} ^{n+1}{(2j-2)!}}.
\end{equation}
The lower half $\odd{R^-}$ can be turned into a determinant in a
manner analogous to Section~\ref{lowsec} (see Figure~\ref{oddpathfi}). 
\begin{figure}
\centertexdraw{
\drawdim cm \linewd.05
\RhombusC 
\move(0 -1)
\RhombusC \RhombusC
\move(0 -2)
\RhombusC \RhombusC \RhombusC
\move(.866025 -.5)
\RhombusB \RhombusA \RhombusB \RhombusA \RhombusB 
\move(1.73205 -1)
\RhombusA \RhombusB \RhombusA \RhombusB 
\move(4.33 -.5)
\RhombusB \RhombusB \RhombusB
\move(3.464 0)
\RhombusC 
\lpatt(.05 .13)
\ringerl(1.299 -.25)
\vdSchritt  \hdSchritt \vdSchritt  \hdSchritt \vdSchritt \fcir f:0 r:.07
\ringerl(2.165063 -.75)
\hdSchritt \vdSchritt  \hdSchritt \vdSchritt \fcir f:0 r:.07 
\ringerl(4.764 -.25)
\vdSchritt \vdSchritt \vdSchritt \fcir f:0 r:.07
}
\caption{\label{oddpathfi} A tiling and lattice paths for $\odd{R^-}$.}
\end{figure}

The starting and end points and the resulting determinant equal
\begin{align*}
\odd{R_i}&=\begin{cases}(2i-1,i+m) &\text {for } i\not= s\\
(2s-2,s+m-1) &\text {for } i=s \end{cases}\\
\odd{S_j}&=(n+j-1,j-1)\\
\odd{A_{ij}}&=\begin{cases}
\dfrac {1} {2} \dbinom{n+m-i}{m+i-j+1} + \dbinom{n+m-i}{m+i-j}
\quad &\text {for } i\not=s\\
\dbinom{n+m-s+1}{m+s-j} &\text {for } i=s.
\end{cases}
\end{align*}
Since the original problem and the claimed final result are invariant
under $n+1-s \rightarrow s$, we can assume $s\not= n$.
(The case $n=s=1$ is trivial to check.)
Then
\begin{align*}
\odd{A_{nj}}&=\frac {1} {2} \binom{m}{m+n-j+1} + \binom{m}{m+n-j}\\
&=\begin{cases}
1\quad \text {for $j=n$,}\\
0\quad \text {else,}
\end{cases}
\end{align*}
since $m$ is a nonnegative integer.
It is easily seen that
$$\odd{A_{ij}}(n,m,s)=A_{ij}(n-1,m+1,s-1)\quad \text {for $i,j< n$,}$$
where $A_{ij}$ is defined in equation~\eqref{aeq}.
We expand $\det_{1\le i,j \le n}(\odd{A_{ij}}(n,m,s))$ with respect
to row $n$ and get
\begin{equation*}
\det_{1\le i,j \le n}(\odd{A_{ij}}(n,m,s))=\det_{1\le i,j \le n-1}
(A_{ij}(n-1,m+1,s-1)),
\end{equation*}
Hence, 
\begin{equation*}
M(\odd{R^-}(n,m,s))=M(R^-(n-1,m+1,s-1)).
\end{equation*}
Thus, Lemma~\ref{g-} yields
\begin{multline} \label{odd-}
M(\odd{R^-})=\frac {2^{\binom{n-2}{2}-1}\h(n-1)(2n-2s-1)!!(2s-3)!!}
{(n-s)! (s-1)! \prod _{i=0} ^{n-2}{(2i+1)!}}\\
\hskip1.5cm\times
\prod _{k=1} ^{n-3}{\left( m+1+k+\frac{1}{2}\right)^{\min(k,n-2-k)} 
\prod_{k=0} ^{n-1}{(m+1+k)^{\min(k+1,n-k)}}}
\end{multline}
Now we substitute the results of
equations~\eqref{odd+} and \eqref{odd-} in \eqref{facto}. We get 
\begin{align}\label{finodd} 
M(\odd{G})&=2^{n-1}M(\odd{R^+})M(\odd{R^-})\\
&=\frac {2^{\binom{n-1}{2}-1}\h(n-1)\h(n+1)(2n-2s-1)!!(2s-3)!!}
{(n-s)! (s-1)! \prod _{j=0} ^{n}{(2j)!}\prod _{i=0}
^{n-2}{(2i+1)!}}\nonumber\\
&\times
\prod _{k=2} ^{n-2}{\left( m+k+\frac{1}{2}\right)^{\min(k-1,n-k-1)} 
\prod_{k=1} ^{n}{(m+k)^{\min(k,n-k+1)}}}
\prod _{2\le i \le j \le n+1} ^{}{(2m+2j-i)},\nonumber
\end{align}
which can easily be transformed to the expression in Theorem~\ref{th2}.
\qed
\end{section}


\end{document}